\documentclass[oneside]{amsart}
\newcommand{\formatswitch}{preprint}

\usepackage{amssymb}
\usepackage{ifthen}
\usepackage{verbatim}
\usepackage{psfig}
\usepackage[all]{xy}

\newcommand{\tref}[1]{(\ref{#1})}

\ifthenelse{\equal{\formatswitch}{big}}{
\setlength{\textwidth}{432pt}
\setlength{\oddsidemargin}{18.8775pt}

\setcounter{tocdepth}{1}
}{
\ifthenelse{\equal{\formatswitch}{paper}}{

\setcounter{tocdepth}{1}
}{
\setcounter{tocdepth}{2}
}
}
\DeclareMathAlphabet\EuScript{U}{eus}{m}{n}
\DeclareMathAlphabet\EuScriptb{U}{eus}{b}{n}

\newcommand{\claimenum}{\renewcommand{\theenumi}{\alph{enumi}}
 \renewcommand{\labelenumi}{\textit{(\theenumi)}}
 \renewcommand{\theenumii}{\roman{enumii}}
 \renewcommand{\labelenumii}{\textit{(\theenumii)}}
 \begin{enumerate}}
\newcommand{\claimenumend}{\end{enumerate}}

\newcommand{\romanenum}{\renewcommand{\theenumi}{\roman{enumi}}
 \renewcommand{\labelenumi}{\textit{(\theenumi)}}
 \renewcommand{\theenumii}{\alph{enumii}}
 \renewcommand{\labelenumii}{\textit{(\theenumii)}}
 \begin{enumerate}}
\newcommand{\romanenumend}{\end{enumerate}}

\newtheorem{dummy}{realdumb}[section]
\newtheorem{thm}{Theorem}

\newtheorem{lemma}[dummy]{Lemma}
\newtheorem{prop}[dummy]{Proposition}
{\theoremstyle{definition} }

{\theoremstyle{definition} }

\renewcommand{\text}{\mathrm}

\newcommand{\strutdepth}{\dp\strutbox}
\newcommand{\marginalnote}[1]
   {\strut\vadjust{\kern-\strutdepth\domarginalnote{#1}}}
\newcommand{\domarginalnote}[1]{\vtop to \strutdepth{
  \baselineskip\strutdepth
   \vss\llap{ #1\ \ }\null}}  
\newcounter{showlabelflag}
\newcounter{makelabelflag}
\newcommand{\showlabels}{\setcounter{showlabelflag}{1}}

\newcommand{\makelabels}{\setcounter{makelabelflag}{1}}
\newcommand{\hidelabels}{\setcounter{showlabelflag}{2}}
\newcommand{\mylabel}[1]{
  \ifthenelse{\value{makelabelflag}=1}
    {\label{#1}}{}
  \ifthenelse{\value{showlabelflag}=1}
    {\marginpar{#1}}{}\relax}

\newcommand{\N}{{\mathbf N}}

\ifthenelse{\equal{\formatswitch}{draft}}{
\showlabels
}{\relax}

\newcommand{\mymargin}[1]{
  \ifthenelse{\value{showlabelflag}=1}
    {\marginpar{#1}}{}\relax}

\newcounter{enumo}

\newcommand{\RRsh}{\kern -1 pt \Rsh}

\newcounter{keepitemnum}

\newcounter{keepitemnumm}

\begin{document}

\bibliographystyle{amsplain}
\begin{center}{\bfseries Presentations of higher dimensional
Thompson groups \footnote{AMS
Classification (2000): primary 20F05, secondary 20M05, 20B27, 
20E32, 20F55, 57S25}}\end{center}
\vspace{3pt}
\begin{center}{MATTHEW G. BRIN}\end{center}
\vspace{4pt}
\vspace{3pt}
\begin{center}July 30, 2003\end{center}

\tableofcontents
\CompileMatrices


\makelabels
\hidelabels


\section{Introduction}

In \cite{brin:hd3} we introduce groups \(nV\) and \(\widehat{nV}\)
for integers \(n\ge1\) and also \(\omega V\) and \(\widehat{\omega
V}\) where \(\omega\) represents the natural numbers.  These are all
subgroups of the homeomorphism group of the Cantor set.  The group
\(1V\) is a group known as the Thompson group \(V\) which is
infinite, simple and finitely presented (see \cite{CFP}).  In
\cite{brin:hd3}, we show that \(2V\) is infinite, simple and
finitely generated and that it is not isomorphic to \(V\), and in
fact not isomorphic to any member of an infinite collection of
infinite, simple, finitely presented groups that are also known as
Thompson groups.

In this paper, we put the group \(2V\) on the same status as \(V\)
and the other infinite, simple, finitely presented Thompson groups
by calculating a finite presentation for \(2V\).  The group \(2V\)
is a subgroup of the (non-simple) group \(\widehat{2V}\) and we also
calculate a finite presentation for \(\widehat{2V}\).  This is done
partly out of necessity since the group \(\widehat{2V}\) is easier
to work with than \(2V\) and a fairly full analysis of
\(\widehat{2V}\) (tantamount to calculating a presentation) must be
done before an analysis of \(2V\) can be done.

The group \(\widehat{2V}\) is a group of fractions of a particularly
nice monoid \(\Pi \) that is also introduced in \cite{brin:hd3}.
The monoid \(\Pi\) can be thought of as a countable, monoid
approximation to the little cubes operads (in dimension 2) of
\cite{may:geom+loop} and \cite{boardman+vogt}.  The significance of
this observation is unclear to the author.  We gain the presentation
and understanding of \(\widehat{2V}\) by first calculating a
presentation for \(\Pi \).

In \cite{brin:hd3}, we give generating sets for \(\Pi \),
\(\widehat{2V}\) and \(2V\).  In that paper, we also list relations
that are satisfied by the elements, but do not show that the
relations suffice to give presentations.  In this paper, we perform
the calculations that show that we have enough relations for a
presentation.

The analysis of the presentation will proceed by using geometric
representations of the elements that distinguish between different
elements, and normal forms for words in the generators that are
derived from the geometric representations.  The work comes in
showing that the relations suffice to reduce an arbitrary word to
one of the normal forms.  One step will make use of the analysis of
\(\widehat{2V}\) by mapping \(\widehat{2V}\) into its subgroup
\(2V\) and using knowlege of \(\widehat{2V}\) to perform part of the
normalization.

The first presentations that we obtain for \(\widehat{2V}\) and
\(2V\) are infinite.  These almost immediately turn into finite
presentations by techniques that have never been formalized into a
single machine, but probably should be.  This machine contains much
of the magical ability of the Thompson groups (with their strong
finiteness properties) to imitate the behavior of much larger
groups.

The monoid \(\Pi \) is not finitely presentable, and in fact is not
finitely generated.

We do not analyze the groups \(nV\) for \(n>2\).  It seems
reasonable to hope that these are also finitely presented.  The
group \(\omega V\) is the ascending union of the groups \(nV\) and
is not finitely generated.

\section{The monoid \protect\(\Pi\protect\)}\mylabel{PiSec}

The monoid \(\Pi\) is defined in \cite{brin:hd3} as a set of
continuous functions from a topological space \(X\) to itself and
also as a set of homeomorphisms from a topological space \(Y\) to
itself.  The invertibility of the homeomorphisms makes it easier to
make the transition to the group of fractions, but the pictures in
the setting of \(X\) are easier to draw.  When discussing \(\Pi\),
we will work with the continous functions on \(X\) and will switch
to the homeomorphisms on \(Y\) when we form the group of fractions.
The shift will be painless.

\subsection{Numbered patterns} The set \(X\) is the union of a
countable set \(\{S_0, S_1, \ldots\}\) of unit squares in the upper
half plane.  The intersection of each \(S_i\) with the \(x\)-axis is
the closed interval \([2i,2i+1]\).

Elements of \(\Pi\) are given by numbered patterns in \(X\).  First
we describe patterns and then we describe numbered patterns.
Patterns are those derivable from a single trivial pattern by a
finite number of simple increments.  We now define these terms.

The \emph{trivial pattern} in \(X\) is the set of rectangles
\(S_i\).  The trivial pattern has exactly one rectangle in each
square \(S_i\).  A \emph{simple increment} to a pattern increases by
one the number of rectangles in a single \(S_i\) by replacing one
rectangle \(R\) in the pattern in \(S_i\) by two congruent
rectangles obtained from \(R\) by dividing \(R\) exactly in half by
either a horizontal line, or a vertical line.  All other rectangles
in the pattern are left alone by the simple increment.  A
\emph{pattern} in \(X\) is a set of rectangles obtainable from the
trivial pattern by a finite number of simple increments.

Note that a pattern in \(X\) (called a \emph{sequence of patterns}
in \cite{brin:hd3}) has more than one rectangle in only finitely
many of the \(S_i\).

A \emph{numbered pattern} is a bijection from
\(\N=\{0,1,2,\ldots\}\) to the rectangles in a pattern for which
there are \(j\) and \(k\) in \(\N\) so that if \(i>k\), then \(S_i\)
has only one rectangle in the pattern and its number is \(i+j\).
The element of \(\N\) associated to a rectangle by the bijection
will be called the \emph{number of the rectangle}.  Below is a
picture of a numbered pattern that is assumed to satisfy the
defining requirement with \(k=3\) and \(j=5\).

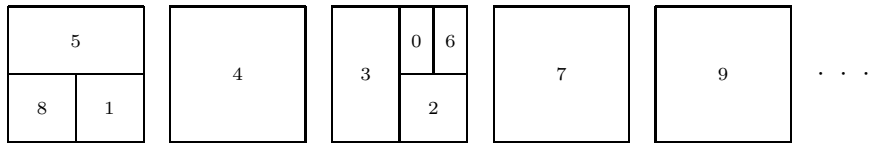
\[
\xy
(-9,-9); (9,-9)**@{-}; (9,9)**@{-}; (-9,9)**@{-}; (-9,-9)**@{-};
(-9,0); (9,0)**@{-};
(0,-9); (0,0)**@{-};
(0,4.5)*{\scriptstyle5};
(-4.5,-4.5)*{\scriptstyle8};
(4.5,-4.5)*{\scriptstyle1};
\endxy
\quad
\xy
(-9,-9); (9,-9)**@{-}; (9,9)**@{-}; (-9,9)**@{-}; (-9,-9)**@{-};
(0,0)*{\scriptstyle4};
\endxy
\quad
\xy
(-9,-9); (9,-9)**@{-}; (9,9)**@{-}; (-9,9)**@{-}; (-9,-9)**@{-};
(0,-9); (0,9)**@{-};
(0,0); (9,0)**@{-};
(4.5,0); (4.5,9)**@{-};
(-4.5,0)*{\scriptstyle3};
(4.5,-4.5)*{\scriptstyle2};
(2.25,4.5)*{\scriptstyle0};
(6.75,4.5)*{\scriptstyle6};
\endxy
\quad
\xy
(-9,-9); (9,-9)**@{-}; (9,9)**@{-}; (-9,9)**@{-}; (-9,-9)**@{-};
(0,0)*{\scriptstyle7};
\endxy
\quad
\xy
(-9,-9); (9,-9)**@{-}; (9,9)**@{-}; (-9,9)**@{-}; (-9,-9)**@{-};
(0,0)*{\scriptstyle9};
\endxy
\quad
\xy
(0,0)*{\cdot};
(3,0)*{\cdot};
(6,0)*{\cdot};
\endxy
\]

\subsection{A monoid of continuous functions}

A numbered pattern determines a continous function \(f\) from \(X\)
to itself which we define separately on each \(S_i\).  If \(R_i\) is
the rectangle in the pattern with number \(i\), then \(f\)
restricted to \(S_i\) is the restriction to \(S_i\) of the unique
affine transformation of the plane of the form \((x,y)\mapsto
(a+2^px, b+2^qy)\) for integers \(p\) and \(q\) that carries \(S_i\)
onto the rectangle \(R_i\).  Note that this carries the lower left
corner of \(S_i\) to the lower left corner of \(R_i\) and so forth.

It is an elementary exercise that the functions corresponding to
numbered patterns in \(X\) form a monoid \(\Pi\) under composition
of functions.  We think of these functions as acting on the left and
we compose from right to left.

Different numbered patterns lead to different functions, so when we
list generators and relations, we will have a criterion for deciding
when two words in the generators give the same element in \(\Pi\).

\subsection{Generators for the monoid}\mylabel{MonGenSec}

The following elements of \(\Pi\) are introduced in \cite{brin:hd3}.
For \(i\ge0\), let \(v_i\) be as pictured below.

\[
\xy
(-6,-6); (-6,6)**@{-}; (6,6)**@{-}; (6,-6)**@{-}; (-6,-6)**@{-};
(0,0)*{\scriptstyle0};
\endxy
\quad
\xy
(-6,-6); (-6,6)**@{-}; (6,6)**@{-}; (6,-6)**@{-}; (-6,-6)**@{-};
(0,0)*{\scriptstyle1};
\endxy
\quad
\xy
(-3,0)*{\cdot};
(0,0)*{\cdot};
(3,0)*{\cdot};
\endxy
\quad
\xy
(-6,-6); (-6,6)**@{-}; (6,6)**@{-}; (6,-6)**@{-}; (-6,-6)**@{-};
(0,0)*{\scriptstyle{i-1}};
\endxy
\quad
\xy
(-6,-6); (-6,6)**@{-}; (6,6)**@{-}; (6,-6)**@{-}; (-6,-6)**@{-};
(0,-6);(0,6)**@{-};
(-3,0)*{\scriptstyle{i}};
(3,0)*{\scriptstyle{i+1}};
\endxy
\quad
\xy
(-6,-6); (-6,6)**@{-}; (6,6)**@{-}; (6,-6)**@{-}; (-6,-6)**@{-};
(0,0)*{\scriptstyle{i+2}};
\endxy
\quad
\xy
(-3,0)*{\cdot};
(0,0)*{\cdot};
(3,0)*{\cdot};
\endxy
\]

In the above picture, each square \(S_j\) with \(j\ne i\) has the
one rectangle, each square \(S_j\) with \(j<i\) is numbered \(j\)
and each square \(S_j\) with \(j>i\) is numbered \(j+1\).

For \(i\ge 0\), let \(h_i\) be as pictured below.

\[
\xy
(-6,-6); (-6,6)**@{-}; (6,6)**@{-}; (6,-6)**@{-}; (-6,-6)**@{-};
(0,0)*{\scriptstyle0};
\endxy
\quad
\xy
(-6,-6); (-6,6)**@{-}; (6,6)**@{-}; (6,-6)**@{-}; (-6,-6)**@{-};
(0,0)*{\scriptstyle1};
\endxy
\quad
\xy
(-3,0)*{\cdot};
(0,0)*{\cdot};
(3,0)*{\cdot};
\endxy
\quad
\xy
(-6,-6); (-6,6)**@{-}; (6,6)**@{-}; (6,-6)**@{-}; (-6,-6)**@{-};
(0,0)*{\scriptstyle{i-1}};
\endxy
\quad
\xy
(-6,-6); (-6,6)**@{-}; (6,6)**@{-}; (6,-6)**@{-}; (-6,-6)**@{-};
(-6,0);(6,0)**@{-};
(0,-3)*{\scriptstyle{i}};
(0,3)*{\scriptstyle{i+1}};
\endxy
\quad
\xy
(-6,-6); (-6,6)**@{-}; (6,6)**@{-}; (6,-6)**@{-}; (-6,-6)**@{-};
(0,0)*{\scriptstyle{i+2}};
\endxy
\quad
\xy
(-3,0)*{\cdot};
(0,0)*{\cdot};
(3,0)*{\cdot};
\endxy
\]

In the above picture, each square \(S_j\) with \(j\ne i\) has the
one rectangle, each square \(S_j\) with \(j<i\) is numbered \(j\)
and each square \(S_j\) with \(j>i\) is numbered \(j+1\).

For \(i\ge 0\), let \(\sigma_i\) be as pictured below.

\[
\xy
(-6,-6); (-6,6)**@{-}; (6,6)**@{-}; (6,-6)**@{-}; (-6,-6)**@{-};
(0,0)*{\scriptstyle0};
\endxy
\quad
\xy
(-6,-6); (-6,6)**@{-}; (6,6)**@{-}; (6,-6)**@{-}; (-6,-6)**@{-};
(0,0)*{\scriptstyle1};
\endxy
\quad
\xy
(-3,0)*{\cdot};
(0,0)*{\cdot};
(3,0)*{\cdot};
\endxy
\quad
\xy
(-6,-6); (-6,6)**@{-}; (6,6)**@{-}; (6,-6)**@{-}; (-6,-6)**@{-};
(0,0)*{\scriptstyle{i-1}};
\endxy
\quad
\xy
(-6,-6); (-6,6)**@{-}; (6,6)**@{-}; (6,-6)**@{-}; (-6,-6)**@{-};
(0,0)*{\scriptstyle{i+1}};
\endxy
\quad
\xy
(-6,-6); (-6,6)**@{-}; (6,6)**@{-}; (6,-6)**@{-}; (-6,-6)**@{-};
(0,0)*{\scriptstyle{i}};
\endxy
\quad
\xy
(-6,-6); (-6,6)**@{-}; (6,6)**@{-}; (6,-6)**@{-}; (-6,-6)**@{-};
(0,0)*{\scriptstyle{i+2}};
\endxy
\quad
\xy
(-3,0)*{\cdot};
(0,0)*{\cdot};
(3,0)*{\cdot};
\endxy
\]

In the above picture, every \(S_j\) has the one rectangle and
every \(S_j\) with \(j\notin\{i,i+1\}\) is numbered \(j\).

The following facts from \cite{brin:hd3} are clear.  We do not
distinguish between a numbered pattern and the element of \(\Pi\)
that it determines.  If \(P\) is a numbered pattern, then \(P_i\) is
the rectangle in \(P\) numbered \(i\).

\begin{lemma}\mylabel{GenPiAction} Let \(P\) be a numbered
pattern.{\claimenum \item The pattern for \(Pv_i\) is gotten from
\(P\) by dividing \(P_i\) vertically, giving the left half the
number \(i\), the right half the number \(i+1\), preserving the
numbers of rectangles \(P_j\) with \(j<i\) and increasing by one the
numbers of rectangles \(P_j\) with \(j>i\).  \item The pattern for
\(Ph_i\) is gotten from \(P\) by dividing \(P_i\) horizontally,
giving the bottom half the number \(i\), the top half the number
\(i+1\), preserving the numbers of rectangles \(P_j\) with \(j<i\)
and increasing by one the numbers of rectangles \(P_j\) with
\(j>i\).  \item The pattern for \(P\sigma_i\) obtained from \(P\) by
exchanging the numbers of rectangles numbered \(i\) and \(i+1\) and
making no other changes.  \item The set \(\{v_i, h_i, \sigma_i \mid
i\in\N\}\) is a  generating set for the monoid
\(\Pi\).  \claimenumend} \end{lemma}

\subsection{Relations for \protect\(\Pi\protect\)}\mylabel{PiRelSec}

The following relations from \cite{brin:hd3} can be checked by hand.
The easiest way is to draw pictures and use Lemma \ref{GenPiAction}.
In \tref{ZSDot} below and in the rest of the paper, we will use the
symbol \(\overline{\sigma}_j\) to refer to the transposition on
\(\N\) that interchanges \(j\) and \(j+1\).

\begin{lemma}\mylabel{PiRels}  The following hold in \(\Pi\).  In
the expressions below, the symbols \(x\) and \(y\) come from \(\{h,v\}\).
\mymargin{PiRelA-F}\begin{alignat}{2} x_jy_i &= y_ix_{j+1},
\qquad&&i<j,\label{PiRelA} \\ \sigma_i^2 &= 1, \qquad &&
i\ge0,\label{PiRelB} \\ \sigma_i\sigma_j &= \sigma_j\sigma_i, \qquad
&& |i-j|\ge2,\label{PiRelC} \\ \sigma_i\sigma_{i+1}\sigma_i &=
\sigma_{i+1}\sigma_i \sigma_{i+1}, \qquad && i\ge0,\label{PiRelD} \\
\sigma_jx_i &= (\sigma_j \cdot x_i)(\sigma_j)^{x_i}, \qquad&&i\ge0,
\,\, j\ge0,\label{PiRelE} \\ v_ih_{i+1}h_i &= h_iv_{i+1}v_i
\sigma_{i+1}, \qquad&& i\ge0,\label{PiRelF} \end{alignat}
where the right side of {\rm \tref{PiRelE}} is given by 
\mymargin{ZSDot}\begin{equation}\label{ZSDot}\sigma_j\cdot x_i =
x_{\overline{\sigma}_j(i)}\end{equation} and 
\mymargin{ZSExp}\begin{equation}
\label{ZSExp}(\sigma_j)^{x_i} = \begin{cases}\sigma_{j+1},
\quad&i<j, \\ \sigma_j\sigma_{j+1}, \quad& i=j, \\
\sigma_{j+1}\sigma_j, \quad &i=j+1, \\ \sigma_j, \quad& i>j+1.
\end{cases}\end{equation}
\end{lemma}

\subsection{First normalization} We keep track of which relations we
use as we simplify a word in the generators of \(\Pi\).  

\begin{lemma}\mylabel{ZSForm} Using the relations \tref{PiRelE}, any
word in the generators from Lemma \ref{GenPiAction}(d) can be
altered to a word in the form \(pq\) where \(p\) is a word in
\(\{v_i, h_i\mid i\in\N\}\) and \(q\) is a word in \(\{\sigma_i\mid
i\in\N\}\).  \end{lemma}

\begin{proof} If we break a word \(w\) in the generators of \(\Pi\)
into a concatenation \(w=p_xp_\sigma r\) where \(p_x\) is the
longest prefix of \(w\) containing nothing but elements of \(\{v_i,
h_i\mid i\in\N\}\) and \(p_\sigma\) is the longest prefix of
\(p_\sigma r\) containing nothing but elements of of
\(\{\sigma_i\mid i\in\N\}\), then we can form a complexity of \(w\)
out of the pair \((a,b)\) where \(a\) is the length of \(p_x\) and
is the most significant part of \((a,b)\), and \(b\) is the length
of \(p_\sigma\).  Complexities are defined to decrease as \(a\)
increases and as \(b\) decreases.  The result follows from the
observation that the relations \tref{PiRelE} do not change the
number of elements of \(\{v_i, h_i\mid i\in\N\}\) that are in a word
and the fact that applications of \tref{PiRelE} lower the
complexity.  \end{proof}

We will need the following trivial strengthening of Lemma
\ref{ZSForm}.

\begin{lemma}\mylabel{RelZSForm} Using the relations \tref{PiRelE},
any word in the generators from Lemma \ref{GenPiAction}(d) of the
form \(pr\) with \(p\) a word in \(\{v_i, h_i\mid i\in\N\}\) can be
altered to a word in the form \(pqs\) where \(q\) is a word in
\(\{v_i, h_i\mid i\in\N\}\) and \(s\) is a word in \(\{\sigma_i\mid
i\in\N\}\).  \end{lemma}

\begin{proof} Apply Lemma \ref{ZSForm} to the subword \(r\).
\end{proof}

\subsection{Labeled, numbered forests}  Elements of \(\Pi\) are
completely classified by numbered patterns, and many words in the
generators of \(\Pi\) lead to one element.  We create a structure
that is intermediate between numbered patterns and words in the
generators from Lemma \ref{GenPiAction}(d).

A forest will be a certain sequence of trees, so we start with
trees.  Our notation is fairly standard and we assume that the
reader is familiar with trees, but we review the terms we will use.

A tree is a non-empty finite set of vertices with two relations
\emph{left child} and \emph{right child}.  Every vertex will either
have one left child and one right child, or it will have no
children.  This makes a tree a binary tree, but all our trees will
be binary and we will not use the adjective ``binary'' when
discussing trees.  A \emph{child} is either a left or right child,
the transitive closure of child is \emph{descendent}, the inverse of
child is \emph{parent} and the transitive closure of parent is
\emph{ancestor}.  Each tree has one vertex, the \emph{root} that is
the ancestor of all other vertices in the tree.

A vertex in a tree is called a \emph{leaf} if it has no children,
and it is called an \emph{interior vertex} otherwise.  The
\emph{trivial tree} has only one vertex which must therefore be both
the root and a leaf.  The trivial tree has no interior vertices.  It
is elementary that the number of leaves of a tree is one more than
the number of interior vertices.

A \emph{labeled tree} is a tree with a label on each interior vertex
where the labels come from \(\{v,h\}\).  A \emph{forest} is a
sequence (indexed over \(\N\)) of trees of which all but finitely
many are trivial.  If we regard the trees of a forest as disjoint,
then we have infinitely many vertices since every tree is non-empty,
but we only have finitely many interior vertices.  A \emph{labeled}
forest is a forest of labeled trees.  The \emph{leaves of the
forest} are elements of the disjoint union of the leaves of the
trees of the forest.  A forest has infinitely many leaves.  If \(F\)
is a forest then \(F_i\) is its \(i\)-th tree.

A \emph{labeled, numbered forest} is a labeled forest \(F\) with a
one-to-one correspondece between \(\N\) and the leaves of the forest
so that there are \(j\) and \(k\) in \(\N\) so that \(i>k\) implies
that \(F_i\) is trivial and its only leaf is numbered \(i+j\).  The
\emph{trivial numbered, labeled forest} is the sequence of trivial
trees so that the sole leaf of the \(i\)-th tree is numbered \(i\).

\subsection{Carets}

In a tree, a triple \((v, v0, v1)\) where \(v\) is an interior
vertex, \(v0\) is the left child of \(v\) and \(v1\) is the right
child of \(v\) is called a \emph{caret}.  The \emph{root} of the
caret \((v, v0, v1)\) is \(v\).  Sending a caret to its root gives a
one-to-one correspondence between the carets in a tree and the
interior vertices of a tree.  We can say that a tree is a union of a
finite number of carets if we are sloppy and declare the trivial
tree to be the union of zero carets.  We talk about labeling
interior vertices, but could just as easily talk about labeling
(roots of) carets.  We introduce carets since they are convenient
when discussing modifications to a tree.

\subsection{Numbered patterns from numbered, labeled
forests}\mylabel{NumPatFromForSec}

Each vertex of a labeled tree corresponds to a unique rectangle in a
unit square.  This can be said inductively.  We
start by declaring that the root corresponds to all of the unit
square.  If a vertex labeled \(v\) corresponds to a rectangle \(R\),
then its left child corresponds to the left half of \(R\) and the
right child corresponds to the right half of \(R\).  If a vertex
labeled \(h\) corresponds to a rectangle \(R\), then its left child
corresponds to the bottom half of \(R\) and the right child
corresponds to the top half of \(R\).

The rectangles corresponding to the leaves of a labeled tree form a
pattern in a unit square.  This is easy to see inductively on the
size of the tree.  A labeled forest thus gives a pattern in \(X\).
A numbered, labeled forest gives a numbered pattern in \(X\) by
giving the rectangle corresponding to a leaf the number of the leaf.
Different numbered, labeled forests can give the same numbered
pattern.  We can discuss examples more easily after the next topic.

\subsection{Numbered, labeled forests from words}\mylabel{BuildForSec}

A word in the generators from Lemma \ref{GenPiAction}(d) determines
a numbered, labeled forest.  We assign the trivial numbered, labeled
forest to the empty word and we define the other assignments
inductively on the length of the word.  If \(w=pa\) with \(a\) from
Lemma \ref{GenPiAction}(d), then the forest \(F\) assigned to \(p\)
is modified depending on \(a\).

If \(a=v_i\), then leaf \(i\) of \(F\) is given two children and the
label \(v\).  We can also describe this as attaching a new caret to
\(F\) by attaching the root of the new caret to leaf \(i\) of \(F\).
The new left child is numbered \(i\), the new right child is
numbered \(i+1\), and each other leaf retains its old number if it
was less than \(i\) and has its number increased by one if the old
number was greater than \(i\).

If \(a=h_i\), then exactly the same thing happens as in the case
\(a=v_i\) except that the new label is an \(h\).

If \(a=\sigma_i\), then the only change to \(F\) is to switch the
numbers of the two leaves that have the numbers \(i\) and \(i+1\).

It is elementary that if \(w\) is a word in the generators from
Lemma \ref{GenPiAction}(d) and \(F\) is the numbered, labeled forest
assigned to \(w\), then the numbered pattern corresponding to \(F\)
is the numbered pattern that determines the same element of \(\Pi\)
as the word \(w\).  There are many words asssigned to the same
numbered, labeled forest.

\begin{lemma}\mylabel{WordForRels} Let \(w\) and \(w'\) be words in
the generators from Lemma \ref{GenPiAction}(d) that are related by
relations \tref{PiRelA}--\tref{PiRelE}.  Then the numbered, labeled
forests assigned to \(w\) and \(w'\) are identical.  \end{lemma}

\begin{proof} Pictures can be drawn for each of the relations.
\end{proof}

The relations \tref{PiRelF} do not preserve the forest.  The
simplest example built from a relation of the form \tref{PiRelF}
shows this and it also gives examples of different numbered, labeled
forests that correspond to the same numbered pattern.

\subsection{Words in the \protect\(v_i\protect\) and
\protect\(h_i\protect\)}

We use the notions of ``confluent'' and ``terminating'' when applied
to relations and rewriting systems.  This material is covered in
numerous places such as \cite{guba+sapir2},
\cite{squier:word-problems} and \cite{brin:zs}.  If we change the
relations in \tref{PiRelA} to
\mymargin{RewriteXY}\begin{equation}\label{RewriteXY}
x_jy_i\longrightarrow y_ix_{j+1}\,\,\,\mathrm{whenever}\,\,\, i<j,
\end{equation} then we have a set of rewriting rules that can be
applied to any word in \(\{v_i, h_i\mid i\in\N\}\) in that we are
allowed to replace a subword like the left side of \tref{RewriteXY}
by the right side of \tref{RewriteXY}, but not the reverse.  It is
an elementary exercise that the rewriting rule \tref{RewriteXY} is
\emph{terminating} in that it cannot be applied an infinite number
of times to a given word, and \emph{locally confluent} in that two
different single applications of \tref{RewriteXY} to a word \(w\) to
give \(w_1\) and \(w_2\) can be ``joined'' by a fourth word \(z\)
that can be obtained from each of \(w_1\) and \(w_2\) by zero or
more applications of \tref{RewriteXY}.  In our situation, getting
\(z\) from \(w_1\) and \(w_2\) will each take no more than two
applications of \tref{RewriteXY}.

The standard fact from such considerations is that the equivalence
class of any word \(w\) in \(\{v_i, h_i\mid i\in\N\}\) under
\tref{PiRelA} contains a unique representative that admits no
applications of \tref{RewriteXY}.  This unique representative is
said to be \emph{irreducible} under \tref{RewriteXY} and gives a
convenient normal form for words in in \(\{v_i, h_i\mid i\in\N\}\).

It is elementary that a word \(x_{i_1}x_{i_2}\cdots x_{i_k}\) where
each \(x\) is chosen independently from \(\{v,h\}\) is irreducible
under \tref{RewriteXY} if and only if \(i_1\le i_2\le \cdots \le
i_k\).  It is now another easy exercise to show that if \(w\) and
\(w'\) are two words in  \(\{v_i, h_i\mid i\in\N\}\) that are
irreducible under \tref{RewriteXY} and are different, then they lead
to different numbered, labeled forests.  Combining this with Lemma
\ref{WordForRels} gives the following.

\begin{lemma}\mylabel{WordsInVHEquiv} Two words in \(\{v_i, h_i\mid
i\in\N\}\) lead to the same labeled, numbered forest if and only if
they are related under \tref{PiRelA}.  \end{lemma}

The leaf numbering of a forest that comes from a word in \(\{v_i,
h_i\mid i\in\N\}\) is particularly simple.  In the next lemma we use
the standard fact that the left-right order on each child pair of a
vertex in a tree leads to a linear, left-right order on the leaves
of a tree.

\begin{lemma}\mylabel{WordsInVHNumber} If the numbered, labeled
forest \(F\) comes from a word in \(\{v_i, h_i\mid i\in\N\}\), then
the leaves of \(F\) are numbered so that the leaves in \(F_i\) have
numbers lower than those in \(F_j\) whenever \(i<j\) and the leaves
in each tree of \(F\) are numbered in increasing order under the
natural left right ordering of the leaves.  \end{lemma}

\subsection{Words leading to the same labeled, numbered forest}

If \(w\) and \(w'\) are words in the generators from Lemma
\ref{GenPiAction}(d) and they are associated to the same labeled,
numbered forest, then we want to conclude that they are related in
some known way.  We know that we can use relations \tref{PiRelE} to
write \(w=pq\) and \(w'=p'q'\) where \(p\) and \(p'\) are words in
\(\{v_i, h_i\mid i\in\N\}\) and \(q\) and \(q'\) are words in
\(\{\sigma_i\mid i\in\N\}\).  Now words in \(\{\sigma_i\mid
i\in\N\}\) can only affect the numbering of a forest, so \(p\) and
\(p'\) must lead to the same labeled forest.  By Lemma
\ref{WordsInVHNumber}, the numbered, labeled forests corresponding
to \(p\) and \(p'\) must be identical.  Thus \(q\) and \(q'\) must
have the same effects on the numbering.  Since the effects of \(q\)
and \(q'\) are permutations calculated from the transpositions
\(\overline{\sigma}_i\), we know that the words \(q\) and \(q'\)
represent the same permutation in the group of finitary permutations
on \(N\).  It is well known that the group of finitary permutations
on \(\N\) are presented by the transpositions
\(\overline{\sigma}_i\) with the relations
\tref{PiRelB}--\tref{PiRelD}.  Thus \(q\) and \(q'\) must be related
by \tref{PiRelB}--\tref{PiRelD}.  Since the numbered, labeled
forests corresponding to \(p\) and \(p'\) are identical, the words
\(p\) and \(p'\) are related by \tref{PiRelA} according to Lemma
\ref{WordsInVHEquiv}.  We have proven the following converse to
Lemma \ref{WordForRels}.

\begin{lemma}\mylabel{TwoWdsSameFor} If two words in the generators
from Lemma \ref{GenPiAction}(d) lead to the same numbered, labeled
forest, then the words are related by \tref{PiRelA}--\tref{PiRelE}.
\end{lemma}

\subsection{Ordering interior vertices}\mylabel{CaretOrderSec}

We will characterize the words in \(\{v_i, h_i\mid i\in \N\}\) that
lead to a given labeled, numbered forest.  Of course, by Lemma
\ref{WordsInVHNumber}, the numbering of such a forest is restricted.

If \(w\) is a word in \(\{v_i, h_i\mid i\in \N\}\) and \(F\) is the
numbered, labeled forest derived from \(w\), then there is a
one-to-one correspondence between the entries in \(w\) and the
interior vertices in \(F\).  The correspondence is easier to
describe by referring to carets instead of interior vertices.  

If \(p_{i-1}\) is the prefix of \(w\) of length \(i-1\) with
\(i\ge1\), and \(F^{i-1}\) is the forest corresponding to
\(p_{i-1}\), then every caret of \(F^{i-1}\) is a caret of \(F^i\)
and every caret of \(F^i\) except one is a caret of \(F^{i-1}\).
Thus it is seen that the set of carets of \(F\) is the ascending
union of the sets of carets of the \(F^i\).  The caret of \(F^i\)
that is not in \(F^{i-1}\) is defined to correspond to the \(i\)-th
entry in \(w\).  Thus the order that the entries appear in \(w\)
gives a linear order to the carets (and thus the interior vertices)
in \(F\).

This linear order respects another order.  It is clear that if
\(w=a_1a_2\cdots a_k\) is a word in \(\{v_i, h_i\mid i\in \N\}\)
with corresponding forest \(F\), then the linear order on the
interior vertices of \(F\) given by the order of the entries in the
word \emph{respects the ancestor relation} in that the interior
vertex for \(a_i\) is never the ancestor of the interior vertex for
\(a_j\) when \(i>j\).  In the proof of the next lemma we treat the
ancestor relation as applied to carets using the one-to-one
correspondence between interior vertices and carets.

\begin{lemma}\mylabel{RealizingIVOrder} If \(F\) is a numbered,
labeled forest with the numbering as in Lemma \ref{WordsInVHNumber},
and if a linear order is given on the interior vertices (and thus of
the carets) of \(F\) that respects the ancestor relation, then there
is a unique word \(w\) in \(\{v_i, h_i\mid i\in \N\}\) leading to
\(F\) so that the order on the interior vertices of \(F\) derived
from the order on the entries in \(w\) is identical to the given
linear order on the interior vertices.  \end{lemma}

\begin{proof} Since the given linear order respects the ancestor
relation, every ancestor of a given caret comes before that caret in
the linear order.  Therefore, there is a sequence of numbered,
labeled forests starting with the trivial forest so that each term
in the sequence is obtained from the previous by adding exactly one
caret and so that the order of addition of carets in this sequence
is exactly the given linear order.  This sequence is unique in that
it is completely determined by the given linear order.  Now a word
can be built up that adds these carets in exactly this order.  Since
each caret with label on its root can be added by exactly one
generator from \(\{v_i, h_i\mid i\in \N\}\), this word is unique.
\end{proof}

\subsection{Secondary labels and normalized
forests}\mylabel{SecondLabelSec}

We have seen that the labels in a labeled forest can be related to
the arrangement of letters in a word.  We will introduce extra
labels to some of the interior vertices in a labeled forest that
correspond to information gathered from the pattern in \(X\) that is
associated to the forest.

Let \(F\) be a labeled forest, and let \(P\) be the pattern
associated to (the leaves of) \(F\).  Numbering will not be
important here.

Let \(u\) be an interior vertex of \(F\).  There is a rectangle
\(R\) that corresponds to \(u\).  Since \(u\) is an interior vertex,
the rectangle \(R\) will be further subdivided by the pattern \(P\).
If the left half of \(R\) is a union of rectangles of \(P\), then it
is necessarily true that the right half of \(R\) is a union of
rectangles of \(P\).  In such case, we say that \(R\) is divided
vertically by \(P\).  Similarly, we say that \(R\) is divided
horizontally by \(P\) if the bottom half of \(R\) is a union of
rectangles of \(P\).  In the case that \(R\) is divided vertically
and also divided horizontally, we say that \(R\) is \emph{fully
divided} and add the \emph{secondary label} ``\(+\)'' to \(u\).
Note that not all interior vertices get secondary labels.

When an interior vertex has a secondary label, then the label from
\(\{v,h\}\) is its \emph{primary label}.  We say that an interior
vertex is \emph{normalized} if it has no secondary label or if its
primary label is \(v\).  We say that a labeled forest \(F\) is
\emph{normalized} if every interior vertex is normalized.

\subsection{Uniqueness of normalized forests}

Uniqueness does not imply existence which will be covered later.
The next lemma establishes uniqueness and is stated so as to be easy
to prove.

\begin{lemma}\mylabel{UniqueNormFor} If two different forests
correspond to the same pattern in \(X\), then at least one of the
two forests is not normalized.  \end{lemma}

\begin{proof} Since different trees in a forest correspond to
patterns in different unit squares in \(X\), we see that it suffices
to look at different labeled trees and assume that they correspond
to the same pattern \(P\) in the unit square.

Pick vertices closest to the root where the two trees differ (as
labeled trees).  A trivial check of cases shows that the vertices
must be interior with different labels.  Since the rectangle
corresponding to a vertex depends only on the labeled path above it
leading to the root and since our choice makes these two labeled
paths the same in the two trees, our differing vertices correspond
to the same rectangle \(R\).  The label must be \(v\) in one tree
and \(h\) in the other, so \(R\) must be fully divided in \(P\) and
the vertices in question must have a secondary label.  Since one of
the trees has \(h\) as the primary label, it is not normalized.
\end{proof}

Existence of a normalized forest for a pattern is a triviality; one
works directly from the pattern.  However, it is not necessary to
argue existence separately since it follows from the next (more
difficult) proposition.

\subsection{Normalized forests from words}\mylabel{NormForSec}

We will prove the following proposition.  If \(w\) is a word in the
generators from Lemma \ref{GenPiAction}(d), then we say the
\emph{length} of \(w\) is the number of appearances in \(w\) of
elements of \(\{v_i, h_i\mid i\in\N\}\).  We note that length is
preserved by the relations \tref{PiRelA}--\tref{PiRelF}.

\begin{prop}\mylabel{NormForFromWord} Let \(w\) be a word in the
generators from Lemma \ref{GenPiAction}(d).  Then \(w\) is related
by \tref{PiRelA}--\tref{PiRelF} to a word corresponding to a
normalized, labeled forest.  \end{prop}

\begin{proof} We will assume that the statement is false for some
word \(w\) of length \(n\) and is true for all words of length less
than \(n\).  There are a number of immediate consequences of this
assumption, not all of which are worth noting.  Two that we need are
that \(n>2\) (since otherwise the corresponding forest has no
secondary labels) and that the following lemma holds.  We will
complete the proof of the proposition after the lemma is stated and
proven.  \end{proof}

\begin{lemma}\mylabel{NormSuffix} Let \(w\) be a word in the
generators from Lemma \ref{GenPiAction}(d) of length \(n\) of the
form \(w=as\) where \(a\) is an element of \(\{v_i, h_i\mid
i\in\N\}\) and \(s\) is of length \(n-1\).  Then we can alter \(w\)
by applications of \tref{PiRelA}--\tref{PiRelF} to \(s\) alone to
give a word whose corresponding forest has all non-root, interior
vertices normalized.  \end{lemma}

\begin{proof} By Lemma \ref{ZSForm}, we can assume that \(s\) is of
the form \(pq\) with \(p\) a word in \(\{v_i, h_i\mid i\in\N\}\) and
\(q\) a word in the \(\sigma_i\) and that \(w=apq\).  Since the
order of the interior vertices of the forest for \(ap\) given by the
order of the letters in \(ap\) must respect the ancestor relation,
we know that the interior vertex corresponding to \(a\) is a root.
By hypothesis, \(s\) is related by \tref{PiRelA}--\tref{PiRelF} to a
word \(s'\) corresponding to a normalized forest.  The pattern \(P\)
for \(as'\) is obtained from the pattern \(P'\) for \(s'\) by
applying the pattern of \(P'\) in unit square \(S_i\) to the
rectangle numbered \(i\) in the pattern for \(a\).  The forest \(F\)
for \(as'\) is obtained from the forest \(F'\) for \(s'\) by
attaching the \(i\)-th tree of \(F'\) to the \(i\)-th leaf of the
forest for \(a\).  Since \(F'\) is normalized, it is seen that \(F\)
has all interior vertices normalized except possibly for the root
vertex of one tree.  \end{proof}

\begin{proof}[Continuation of the proof of Proposition
\ref{NormForFromWord}] By Lemma \ref{NormSuffix}, we can assume that
the forest \(F\) for \(w\) has all non-root interior vertices
normalized.  Let \(P\) be the pattern corresponding to \(w\).

By Lemmas \ref{ZSForm} and \ref{WordForRels}, we can assume that
\(w\) is of the form \(w=pq\) with \(p\) a word in \(\{v_i, h_i\mid
i\in\N\}\) and \(q\) a word in the \(\sigma_i\).  Let \(u\) be a
root of \(F\) that is not normalized and let \(u_0\) and \(u_1\) be
the left and right children, respectively, of \(u\).  By Lemma
\ref{RealizingIVOrder}, we can assume that the first three letters
of \(p\) correspond, in order, to \(u\), \(u_1\) and \(u_0\).  This
choice of order is deliberate.

Since \(u\) is not normalized, its corresponding rectangle \(R\) is
fully divided and its label is \(h\).  Thus the rectangles
corresponding to \(u_0\) and \(u_1\) are the bottom and top
rectangles, respectively, of \(R\) and must be both vertically
divided since \(R\) is fully divided.  Since \(u_0\) and \(u_1\) are
normalized, a quick check of cases shows that they are both labeled
\(v\).  Thus the first three letters of \(p\) are \(h_iv_{i+1}v_i\).
Using \tref{PiRelF}, these can be replaced by
\(v_ih_{i+1}h_i\sigma_{i+1}\).  This gives a word \(w'\) whose
corresponding forest \(F'\) is different from \(F\) but whose
corresponding pattern is still \(P\).

The root corresponding to the initial letter \(v_i\) is now
normalized.  Now a second application of Lemma \ref{NormSuffix}
normalizes all other vertices.  The result is a normalized forest.
\end{proof}

\subsection{A presentation for \protect\(\Pi\protect\)}  We are
ready for the following.

\begin{thm}\mylabel{PiPresentation} The monoid \(\Pi\) is presented
by using the generators from Lemma \ref{GenPiAction}(d) and
relations \tref{PiRelA}--\tref{PiRelF}.  \end{thm}

\begin{proof} Let two words give the same element of \(\Pi\).  Since
elements correspond to numbered patterns, they give the same
numbered pattern.  From Proposition \ref{NormForFromWord} and Lemma
\ref{UniqueNormFor}, we can assume that the two words correspond to
the same labeled forest.  Since the forests correspond to the same
numbered pattern, the forests must have the same numbering.  The
result now follows from Lemma \ref{TwoWdsSameFor}.  \end{proof}

\section{The group \protect\(\widehat{2V}\protect\)}

The group \(\widehat{2V}\) is the group of right fractions of
\(\Pi\).  However, as in \cite{brin:hd3} it is easier to change the
representation of \(\Pi\) to make the elements invertible than it is
to apply the usual theorem (Ore's theorem, Theorem 1.23 of
\cite{cliff+prest:I}).

\subsection{New patterns} To invert the elements of \(\Pi\), we
alter the meaning of vertical and horizontal divisions slightly.  If
\(R\) is a rectangle, than the new notion of vertical division
replaces \(R\) by its left third and right third.  The new
horizontal division of \(R\) replaces \(R\) by its bottom third and
top third.  Patterns defined with these steps do not give
collections of rectangles that cover all of the unit squares in
\(X\).  However, a pattern defined this way will cover copies of
\(C\times C\) in \(X\) where \(C\) is the Cantor set defined in the
usual way as the ``deleted middle thirds'' set in the unit interval
and this inclusion of \(C\) in the unit interval \(I\) induces the
natural inclusion of \(C\times C\) in the unit square \(I\times I\).
The covering of \(C\times C\) will be by pairwise disjoint closed
and open sets in \(C\times C\).  Numberings will be handled in the
same way in the new and old patterns and will be made to correspond.

The following is a picture from \cite{brin:hd3} that shows how an
old numbered pattern in a unit square converts to a new numbered
pattern covering \(C\times C\) in that square.

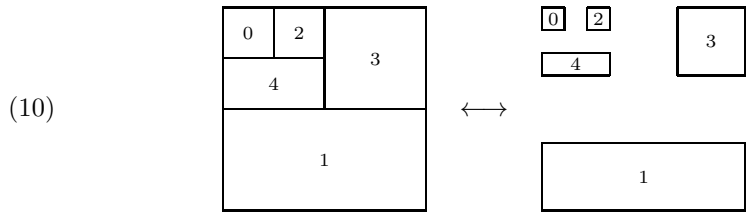
\mymargin{OldToNewPat}\begin{equation}\label{OldToNewPat}
\xy (-13.5,-13.5); (-13.5,13.5)**@{-}; (13.5,13.5)**@{-}; (13.5,-13.5)**@{-};
(-13.5,-13.5)**@{-};
(-13.5,0); (13.5,0)**@{-}; (0,0); (0,13.5)**@{-};
(-13.5,6.75); (0,6.75)**@{-};
(-6.75,6.75); (-6.75,13.5)**@{-};
(0,-6.75)*{\scriptstyle1};
(6.75,6.75)*{\scriptstyle3};
(-6.75,3.375)*{\scriptstyle4};
(-10.125,10,125)*{\scriptstyle0};
(-3.375,10,125)*{\scriptstyle2};
\endxy
\quad
\longleftrightarrow
\quad
\xy
(-13.5,-13.5); (-13.5,-4.5)**@{-}; (13.5,-4.5)**@{-}; 
(13.5,-13.5)**@{-}; (-13.5,-13.5)**@{-};
(0,-9)*{\scriptstyle1};
(-13.5,4.5); (-13.5,7.5)**@{-}; (-4.5,7.5)**@{-}; 
(-4.5,4.5)**@{-}; (-13.5,4.5)**@{-};
(-9,6)*{\scriptstyle4};
(-13.5,10.5); (-13.5,13.5)**@{-}; (-10.5,13.5)**@{-}; 
(-10.5,10.5)**@{-}; (-13.5,10.5)**@{-};
(-12,12)*{\scriptstyle0};
(-7.5,10.5); (-7.5,13.5)**@{-}; (-4.5,13.5)**@{-}; 
(-4.5,10.5)**@{-}; (-7.5,10.5)**@{-};
(-6,12)*{\scriptstyle2};
(13.5,4.5); (13.5,13.5)**@{-}; (4.5,13.5)**@{-}; 
(4.5,4.5)**@{-}; (13.5,4.5)**@{-};
(9,9)*{\scriptstyle3};
\endxy
\end{equation}

\subsection{Elements of the group} We put a copy of \(C\times C\) in
each unit square \(S_i\) of \(X\) and let \(Y\) be the union of
these copies of \(C\times C\).  Elements of \(\widehat{2V}\) are
self homeomorphisms of \(Y\) and are defined by pairs of numbered
patterns \((P,Q)\).  We think of \(P\) as the range pattern and
\(Q\) as the domain pattern.  We adopt this convention to make
formulas for composition look nicer.  For each \(i\) in \(\N\), the
homeomorphism defined by \((P,Q)\) takes the intersection of \(Y\)
with the \(i\)-th rectangle under \(Q\) onto the intersection \(Y\)
with the \(i\)-th rectangle of \(P\) by the unique affine
transformation \((x,y)\mapsto (a+3^jx\,,\, b+3^ky)\) with \(j\) and
\(k\) integers that does so.

It turns out that many pairs will represent the same element of
\(\widehat{2V}\).  For example, the identity of \(\widehat{2V}\) is
represented by all pairs of the form \((P,P)\).

In spite of the fact that the closed and open sets covering
\(C\times C\) look like the right side of \tref{OldToNewPat}, we
will continue to think of patterns as drawn in the left side of
\tref{OldToNewPat}.  Thus we continue to talk about rectangles being
divided in half and not thirds.  When a pattern based on halves is
used to create an element of \(\widehat{2V}\), it has to be
converted first into a pattern based on thirds.

\subsection{Group of right fractions} It is more useful to view
elements as compositions of two homeomorphisms.  If \(E\) is the
trivial pattern, then \((P,Q)\) is the composition of \((P,E)\) and
\((E,Q)\).

With the homeomorphisms acting on the left, we write \[(P,Q)=
(P,E)(E,Q) = (P,E)(Q,E)^{-1}.\] The elements of the form \((P,E)\)
create a copy of the monoid \(\Pi\).  Specifically, sending \(P\)
interpreted as a pattern in \(X\) to the pair \((P,E)\) in
\(\widehat{2V}\) creates an isomorphic embedding of \(\Pi\) into
\(\widehat{2V}\).  If we identify the element \(P\) of \(\Pi\) with
\((P,E)\) in \(\widehat{2V}\), then this establishes \(\Pi\) as a
monoid in \(\widehat{2V}\) with the property that every element of
\(\widehat{2V}\) is of the form \(PQ^{-1}\) with both \(P\) and
\(Q\) in \(\Pi\).  From \cite{cliff+prest:I} (Page 36 and Problem 3
of Page 37), this establishes \(\widehat{2V}\) as a group of
right fractions of \(\Pi\).

From this point, we will use both \((P,Q)\) and \(PQ^{-1}\) to
denote the same element of \(\widehat{2V}\).  Note that if \(M\) is
any other element of \(\Pi\), then \((PM,QM)\) represents the same
element \((PM)(QM)^{-1}=PQ^{-1}\) as \((P,Q)\).

\subsection{Presentation of \protect\(\widehat{2V}\protect\)}  By a
well known extension of Ore's theorem (see, for example, Proposition
2.4 of \cite{brin:bv}), a monoid presentation for \(\Pi\) is a group
presentation for \(\widehat{2V}\).  Thus we have the following.

\begin{thm}\mylabel{WidehatTwoVPres} The group \(\widehat{2V}\) is
presented by using the generators from Lemma \ref{GenPiAction}(d)
and relations \tref{PiRelA}--\tref{PiRelF}.  \end{thm}

\subsection{A semi-normal form}

There will be a natural homomorphism from \(\widehat{2V}\) into
\(2V\) (actually and embedding, but the injective property will not
be needed).  We will use this homomorphism and the following lemma
to prove consequences about the relations in \(2V\).  The lemma is
stated in a way that will be easy to use.

\begin{lemma}\mylabel{WidehatTwoVForm} Let \(w\) be a word in the
generators from Lemma \ref{GenPiAction}(d) and their inverses.
By applying the relations \tref{PiRelA}--\tref{PiRelF} to \(w\), we
can obtain a word of the form \(LMR\) where \(L\) and \(R^{-1}\) are
words in \(\{v_i, h_i\mid i \in \N\}\) and \(M\) is a word in
\(\{\sigma_i\mid i \in \N\}\).  \end{lemma}

\begin{proof} Since the relations \tref{PiRelA}--\tref{PiRelF} are
all of the presenting relations of \(\widehat{2V}\), we only have to
argue that every element of \(\widehat{2V}\) has a representative in
the desired form.  But every element of \(\widehat{2V}\) can be put
in the form \(pq^{-1}\) with \(p\) and \(q\) a word in the
generators from Lemma \ref{GenPiAction}(d) and Lemma \ref{ZSForm}
puts each of \(p\) and \(q\) in the form \(ab\) with \(a\) a word in
\(\{v_i, h_i\mid i \in \N\}\) and \(b\) a word in \(\{\sigma_i\mid i
\in \N\}\).  The fact that elements of \(\{\sigma_i\mid i \in \N\}\)
are their own inverses completes the proof.  \end{proof}

\subsection{An interchange formula}  Unfortunately we will need more
detail than supplied by Lemma \ref{WidehatTwoVForm}.  The next
lemma allows us to predict to some extent what we might see when
putting a word in the generators from Lemma \ref{GenPiAction}(d)
into semi-normal form.

\begin{lemma}\mylabel{MoveALeftNegWord} Let \(w\) be a word in
\(\{\sigma_i, v_i^{-1}, h_i^{-1}\mid i\in\N\}\).  Then applications
of \tref{PiRelA}--\tref{PiRelF} can be use to put \(wv_i\) in the
form \(pw'\) where \(p\) is a word in \(\{v_i\mid i \in\N\}\) and
\(w'\) is a word in \(\{\sigma_i, v_i^{-1}, h_i^{-1}\mid
i\in\N\}\).  \end{lemma}

\begin{proof} Recall that elements of \(\widehat{2V}\) of the form
\((P,E)\), where \(E\) is the trivial pattern, form a copy of the
monoid \(\Pi\).  The word \(w\) is the inverse of a word in this
copy of \(\Pi\) and is of the form \((E,P)\).  The element \(v_i\)
is of the form \((V_i,E)\) where we use \(V_i\) to denote the first
pattern shown in Section \ref{MonGenSec}.

We can multiply the elements \(w=(E,P)\) and \(v_i=(V_i,E)\) if we
can get a pairs representing \(w\) and \(v_i\) so that the second
pattern for \(w\) is the same as the first pattern for \(v_i\).  As
pairs in a group of fractions, we can get new pairs from old by
multiplying both entries on the right by the same thing.  We will
use Lemma \ref{GenPiAction} to understand this right multiplication.

We start to get a pattern that is a common right multiple of \(P\)
and \(V_i\) by superposing the two patterns.  Since \(V_i\) has only
one non-trivial division consisting of a vertical line in the
\(i\)-th square, we need only draw a vertical line in the \(i\)-th
square of the pattern \(P\).  This might have the consequence of
vertically dividing several of the rectangles of \(P\).  According
to Lemma \ref{GenPiAction}, this can be accomplished by multiplying
the element of \(\Pi\) given by \(P\) on the right by a word \(p\)
in the \(v_i\).  Abusing notation somewhat, we now use the pair
\((Ep,Pp)\) to represent \(w\).  Since \(E\) is the trivial pattern,
we get \(w=(p,Pp)\).

Now \(V_i\) must be converted to \(Pp\).  This is accomplished by
subdividing various rectangles in \(V_i\) and applying a word in the
\(\sigma_i\) to get the right numbering.  According to Lemma
\ref{GenPiAction}, a word \(w'\) in \(\{\sigma_i, v_i, h_i \mid i
\in \N\}\) must be applied to the right of \(V_i\).  Thus we get
\(V_iw'=Pp\) and \(v_i=(V_iw',Ew')=(Pp,w')\) and
\[wv_i=(p,Pp)(Pp,w') = (p,w') = (p,E)(E,w').\] This is exactly what
was wanted.  \end{proof}


\section{The group \protect\(2V\protect\)}

\subsection{The elements} Recall that \(X\) is the disjoint union of
the squares \(S_i\), \(i\in\N\), and that \(Y\) is a subset of
\(X\).  The group \(2V\) is the subgroup of \(\widehat{2V}\)
consisting of those elements that act as the identity off \(Y\cap
S_0\).  It is easy to see that these are the elements that are
representable by pairs \((P,Q)\) for which there is an \(n\in \N\)
so that each of \(P\) and \(Q\) satisfy the following: (1) the
number of rectangles in \(S_0\) is \(n\), (2) the number of
rectangles in each \(S_i\), \(i>0\), is one, and (3) the number of
the rectangle in \(S_i\) is \(i+n-1\) for \(i>0\).

\newcommand{\opi}{\overline{\pi}}

\subsection{The generators}\mylabel{TwoVGenSec}

The following elements of \(2V\) are shown
in \cite{brin:hd3} to generate \(2V\).  
\begin{alignat*}{2}
A_i &=(v_0^{i+1}v_1, \, v_0^{i+2}), &\qquad&i\ge0, \\
B_i &=(v_0^{i+1}h_1, \, v_0^{i+2}), &&i\ge0, \\
C_i &=(v_0^ih_0,\, v_0^{i+1}), &&i\ge0, \\
\pi_i &= (v_0^{i+2}\sigma_1, \, v_0^{i+2}), && i\ge0, \\
\opi_i &= (v_0^{i+1}\sigma_0, \, v_0^{i+1}), &&i\ge0.
\end{alignat*}  We let \[\Sigma = \{A_i, B_i, C_i, \pi_i, \opi_i\mid
i\in \N\}.\]

The argument that \(\Sigma\) is a generating set for \(2V\) is not
relevant to this paper.

\subsection{The relations}\mylabel{TwoVRelSec}

The list of relations is longer.  In \cite{brin:hd3} it is argued
that the following relations hold in \(2V\) where 
\(X\) and \(Y\) represent symbols from \(\{A,B\}\).
\mymargin{TwoVRelA-Q}\begin{alignat}{2} 
X_qY_m &= Y_mX_{q+1},&\quad&m<q, \label{TwoVRelA}\\
\pi_qX_m &= X_m\pi_{q+1}, &&m<q, \label{TwoVRelB}\\
\pi_qX_q &= X_{q+1}\pi_q\pi_{q+1}, &&q\ge0, \label{TwoVRelC}\\
\pi_qX_m &= X_m\pi_q, && m>q+1, \label{TwoVRelD}\\
\opi_qX_m &= X_m\opi_{q+1}, &&m<q, \label{TwoVRelE}\\
\opi_mA_m &= \pi_m\opi_{m+1}, &&m\ge0, \label{TwoVRelF}\\
\opi_mB_m &= C_{m+1} \pi_m \opi_{m+1}, &&m\ge0, \label{TwoVRelG}\\
C_qX_m &= X_mC_{q+1}, &&m<q, \label{TwoVRelH}\\
C_mA_m &= B_mC_{m+2}\pi_{m+1}, &&m\ge0, \label{TwoVRelI}\\
\pi_qC_m &= C_m\pi_q, && m>q+1, \label{TwoVRelJ}\\
A_mB_{m+1}B_m &= B_mA_{m+1}A_m \pi_{m+1}, &&m\ge0, \label{TwoVRelK}\\
\pi_q\pi_m &= \pi_m\pi_q, 
  &&|m-q|\ge2, \label{TwoVRelL}\\
\pi_m\pi_{m+1}\pi_m &= \pi_{m+1}\pi_m\pi_{m+1},
  &&m\ge0, \label{TwoVRelM}\\
\opi_q\pi_m &= \pi_m\opi_q, 
  &&q\ge m+2, \label{TwoVRelN}\\
\pi_m\opi_{m+1}\pi_m &= \opi_{m+1}\pi_m\opi_{m+1},
  &&m\ge0, \label{TwoVRelO}\\
\pi_m^2 &= 1, &&m\ge0, \label{TwoVRelP}\\
\opi_m^2 &= 1, &&m\ge0.\label{TwoVRelQ}
\end{alignat}
It is our task to show that the generators in \(\Sigma\) and the
relations \tref{TwoVRelA}--\tref{TwoVRelQ} present \(2V\).

We could eliminate the generators \(C_i\) by using \[C_m =
(\opi_mB_m\opi_{m+1}\pi_m)(B_m\pi_{m+1}A_m^{-1}),\] as is shown in
\cite{brin:hd3}.  However this does not seem to simplify the
calculations below.

\subsection{Strategy}

We need only show that if a word represents the trivial element,
then the word is reducible to the trivial word by the relations 
\tref{TwoVRelA}--\tref{TwoVRelQ}.  However, it is hard to use the
fact that the element represented is trivial until the word has been
simplified significantly.  Thus we reduce an arbitrary word to a
particularly nice form first, and then take into account that the
represented element is trivial.

\subsection{Conventions}\mylabel{ConventionSec}

We let \(G\) be the group with generators
from Lemma \ref{GenPiAction}(d) and relations
\tref{TwoVRelA}--\tref{TwoVRelQ}.  If two words \(w\) and \(w'\) in
the generators of \(G\) represent the same element of \(G\), then we
will write \(w\sim w'\).  

We will be giving different treatment to the positive and negative
powers of the generators.  Thus from now on we will work with the
generating set \[\Sigma_s = \{A_i, B_i, C_i, \pi_i, \opi_i,
A_i^{-1}, B_i^{-1}, C_i^{-1}\mid i\in\N\}\] and treat it as a group
of semigroup generators of the group \(G\).  We will have no need to
distinguish between \(\pi_i\) and \(\pi_i^{-1}\) or between
\(\opi_i\) and \(\opi_i^{\,-1}\) because of the relations
\(\pi_i^2=1\) and \(\opi_i^{\,2}=1\).

We will never have reason to discuss subsets of \(\Sigma_s\) on the
basis of the values of the subscripts.  Thus we will often refer to
subsets of \(\Sigma_s\) by leaving out the subscripts and referring
to words in these subsets in the following form.  If \(S\) is a
subset of the symbols \[\{A, B, C, \pi, \opi, A^{-1}, B^{-1},
C^{-1}\},\] then we will write \(w(S)\) to indicate a word in the
symbols from \(S\), subscripted with values from \(\N\).  For
example, \(w(A, \pi, B^{-1})\) refers to a word in the subset
\(\{A_i, \pi_i, B^{-1}_i\mid i \in\N\}\) of \(\Sigma_s\).

\subsection{The \protect\(LMR\protect\) form, Part I} 

Let \(w\) be a word in \(\Sigma_s\).  Our first task will be to show
that \(w\sim LMR\) where \(L\) and \(R^{-1}\) are words of the form
\(w(A, B, C)\) and \(M\) is a word of the form \(w(\pi, \opi)\).
The calculations that do this are rather intricate and will be done
in several steps.  We will start with words in a specific subset of
the generators.  Then we will add generators one type at a time.
The initial argument will be based on what we know about
\(\widehat{2V}\).  The remaining arguments will be detailed
calculations based on the relations
\tref{TwoVRelA}--\tref{TwoVRelQ}.

\begin{lemma}\mylabel{LMRABpi} (1) Let \(w\) be of the form
\(w(A,B,\pi,A^{-1}, B^{-1})\).  Then \(w\sim LMR\) where \(L\) and
\(R^{-1}\) are words of the form \(w(A,B)\) and \(M\) is of the form
\(w(\pi)\).

(2) Let \(w\) be of the form \(w(A,B,\pi)\).  Then
\(w\sim LM\) where \(L\) is a word of the form
\(w(A,B)\) and \(M\) is of the form \(w(\pi)\).  \end{lemma}

\begin{proof} There is a homomorphism from \(\widehat{2V}\) to \(G\)
defined by \(v_i\mapsto A_i\), \(h_i\mapsto B_i\), \(\sigma_i\mapsto
\pi_i\).  This is seen since the relations of \(\widehat{2V}\)
correspond to the relations of \(G\) according to the following
table \[\begin{split} \tref{PiRelA}&\rightarrow \tref{TwoVRelA}, \\
\tref{PiRelB} &\rightarrow \tref{TwoVRelP}, \\ \tref{PiRelC}
&\rightarrow \tref{TwoVRelL}, \end{split}\qquad\qquad\qquad
\begin{split} \tref{PiRelD} &\rightarrow \tref{TwoVRelM}, \\
\tref{PiRelE} &\rightarrow \tref{TwoVRelB}\hbox{--}\tref{TwoVRelD},
\\ \tref{PiRelF} &\rightarrow \tref{TwoVRelK}.  \end{split}\] In the
correspondence \tref{PiRelE} with \tref{TwoVRelB}--\tref{TwoVRelD},
we make use of the fact that \(\pi_qX_q = X_{q+1}\pi_q\pi_{q+1}\)
implies \(\pi_qX_{q+1} = X_q\pi_{q+1}\pi_q\) since the \(\pi_i\) are
their own inverses in \(G\).

The statement (1) now follows from Lemma \ref{WidehatTwoVForm}.  The
statement (2) follows from Lemma \ref{ZSForm} and the corresponding
monoid homomorphism from \(\Pi\) to \(G\).
\end{proof}

We note that a direct proof of Lemma \ref{LMRABpi}(1) from the
relations seems rather complicated.

\subsection{Subscript raising formulas}

Because of the dependence of some relations on relative values of
subscripts, it will be convenient to alter some subscripts.  The
next lemma allows a subscripted generator to be replaced by the same
generator with a higher subscript at the expense of introducing
words in the other generators.

\begin{lemma}\mylabel{Raises}  The following are consequences of the
relations \tref{TwoVRelA}--\tref{TwoVRelQ},
\[\begin{split}
  C_r &\sim C_{r+1}B_r \pi_{r+1} A_r^{-1} \\
  \opi_r &\sim \pi_r \opi_{r+1} A_r^{-1} \\ 
     &\sim A_r \opi_{r+1}\pi_r
  \end{split}
\]  \end{lemma}

\begin{proof} The first follows from \(C_mA_m\sim
B_mC_{m+2}\pi_{m+1}\) and \(C_qB_m\sim B_mC_{q+1}\) when \(m<q\).
The second and third follow from \(\opi_mA_m\sim \pi_m\opi_{m+1}\)
and \(\opi_m^{\,2}\sim 1\).  \end{proof}

\subsection{Interchanges}

The basic tools for getting words into nicer form will be
``reversals'' of generators that are in the wrong order.  If \(LMR\)
form is desired, then the apperance of \(\pi_q A_r\) in a word
will an obstruction to getting this form.   The resolution will
depend on the relative values of \(q\) and \(r\).  For example if
\(r<q\), then we can replace the letters with \(A_r\pi_{q+1}\).
However, if \(r=q\), then we get \(A_{q+1}\pi_q\pi_{q+1}\).

As can be seen, sometimes an interchange results in a word that is
fairly complex.  The above examples are quite simple and the
interchanges get considerably worse.  It is often more important to
know the form that results from an interchange than the actual value
of the word.  Thus for example, we can write \(\pi_qC_r \sim C
w(A^{-1}, \pi, B)\) when \(r<q+2\), rather than the more exact and
complicated \[\pi_qC_r \sim C_{q+2}\pi_q
(B_{q+1}\pi_{q+2}A_{q+1}^{-1})(B_{q}\pi_{q+1}A_{q}^{-1}) \cdots
(B_{r}\pi_{r+1}A_{r}^{-1}).\] The omission of the subscript of \(C\)
on the right side of \(C\) in \(\pi_qC_r \sim C w(A^{-1}, \pi, B)\)
is deliberate since its exact value will not be important.

\subsection{Interchange formulas}

Below we give the formulas that we need to get a word into \(LMR\)
form. 

Our notation is best illustrated by example.  In writing
\[B_q^{-1}A_r \sim \begin{cases} AB^{-1}, & r\ne q, \\
w(A) \pi w(B^{-1}), &r=q.\end{cases}\] we say that the expression on
the left can be replaced by the expressions on the right under the
conditions stated.  The subscripts are on the right are left
unspecified as they will not be important.

We separate the formulas for moving the different generators to make
them eaiser to refer to.

\begin{lemma}\mylabel{MoveALeft} The following are consequences of
the relations \tref{TwoVRelA}--\tref{TwoVRelQ} and are used to move
positive powers of \(A\) to the left.  Their inverses can be used to
move negative powers of \(A\) to the right.  In the last formula,
the two words of form \(w(\pi,A^{-1},B^{-1})\) are not to be assumed
identical.
\[
\begin{split}
A_q^{-1}A_r &\sim 
  \begin{cases} 
    AA^{-1}, & r\ne q, \\ 
    1, &r=q.
  \end{cases} \\ 
B_q^{-1}A_r &\sim 
  \begin{cases} 
    AB^{-1}, &r\ne q, \\ 
    w(A)\pi w(B^{-1}), &r=q. 
  \end{cases} \\ 
C_q^{-1}A_r &\sim
  \begin{cases} 
    AC^{-1}, & r<q, \\ 
    w(A, \pi, B^{-1})C^{-1}, &r\ge q. 
  \end{cases} \\ 
\pi_q A_r & \sim A w(\pi). \\ 
\opi_q A_r &\sim
  \begin{cases} 
    A\opi, &r<q, \\ 
    \pi\opi, &r=q, \\ 
    w(A) \opi w(\pi), &r>q. 
  \end{cases} \\
w(\pi,A^{-1},B^{-1})A_r &\sim
  w(A)w(\pi,A^{-1},B^{-1}).
\end{split}
\] 
\end{lemma}

\begin{proof} The last formula follows from Lemma
\ref{MoveALeftNegWord} exactly as Lemma \ref{LMRABpi} follows from
Lemma \ref{WidehatTwoVForm}.  For the rest, we will discuss the less
simple instances and leave the others to the reader.

For \(C_q^{-1}A_r\) with \(r\ge q\), we use the inverse of the first
line in Lemma \ref{Raises} repeatedly to get
\[
\begin{split}
C_q^{-1}A_r 
&=
(A_q\pi_{q+1}B_q^{-1})(A_{q+1}\pi_{q+2}B_{q+1}^{-1})
\cdots (A_r\pi_{r+1}B_r^{-1})C_{r+1}^{-1}A_r
\\ &=
(A_q\pi_{q+1}B_q^{-1})(A_{q+1}\pi_{q+2}B_{q+1}^{-1})
\cdots (A_r\pi_{r+1}B_r^{-1})A_rC_{r+2}^{-1}.
\end{split}
\]

For \(\opi_qA_r\) with \(r>q\), we use the third line in Lemma
\ref{Raises} repeatedly to write
\[
\begin{split}
\opi_qA_r 
&= 
A_qA_{q+1}\cdots A_{r-1}\opi_r \pi_{r-1}\pi_{r-2}\cdots \pi_q A_r
\\ &= 
A_qA_{q+1}\cdots A_{r-1}\opi_r \pi_{r-1} A_r \pi_{r-2}
\pi_{r-3}\cdots \pi_q
\\ &= 
A_qA_{q+1}\cdots A_{r-1}\opi_r A_{r-1} \pi_r \pi_{r-1}
\pi_{r-2} \pi_{r-3}\cdots \pi_q
\\ &= 
A_qA_{q+1}\cdots A_{r-1}A_{r-1} \opi_{r+1}
\pi_r \pi_{r-1} \cdots \pi_q
\end{split}
\]
\end{proof}

\begin{lemma}\mylabel{MoveBLeft} The following are consequences of
the relations \tref{TwoVRelA}--\tref{TwoVRelQ} and are used to move
positive powers of \(B\) to the left.  Their inverses can be used to
move negative powers of \(B\) to the right.  \[\begin{split}
A_q^{-1}B_r &\sim \begin{cases} BA^{-1}, &r\ne q, \\ w(B)\pi
w(A^{-1}), &r=q. \end{cases} \\ B_q^{-1}B_r &\sim \begin{cases}
BB^{-1}, & r\ne q, \\ 1, &r=q.\end{cases} \\ C_q^{-1}B_r &\sim
\begin{cases} BC^{-1}, & r<q, \\ w(A, \pi, B^{-1})C^{-1}, &r\ge
q. \end{cases} \\ \pi_q B_r & \sim B w(\pi). \\ \opi_q B_r &\sim
\begin{cases} B\opi, &r<q, \\ C \pi\opi, &r=q, \\ w(A) B \opi
w(\pi), &r>q. \end{cases} \end{split}\] \end{lemma}

\begin{proof} The proof differs little from that of Lemma
\ref{MoveALeft}  \end{proof}

\begin{lemma}\mylabel{MoveCLeft} The following are consequences of
the
relations \tref{TwoVRelA}--\tref{TwoVRelQ} and are used to move
positive powers of 
\(C\) to the left.  Their inverses can be used to move negative
powers of \(C\) to the right.  
\[
\begin{split}
A_q^{-1}C_r
&\sim
\begin{cases} CA^{-1}, &q<r, \\ Cw(A^{-1}, \pi, B), &q\ge r.
\end{cases} \\
B_q^{-1}C_r
&\sim
\begin{cases} CB^{-1}, &q<r, \\ Cw(A^{-1}, \pi, B), &q\ge r.
\end{cases} \\
C^{-1}_qC_r 
&\sim 
\begin{cases} 
w(A^{-1}, \pi, B),   &r<q, \\
1                    &r=q, \\
w(A, \pi, B^{-1}),   &r>q.
\end{cases}  \\
\pi_qC_r
&\sim
\begin{cases} C\pi, & r>q+1, \\ Cw(A^{-1}, \pi, B), &r\le q+1.
\end{cases} \\
\opi_qC_r
&\sim
\begin{cases} B\opi\pi, &r=q+1, \\
              w(A)B\opi w(\pi), &r>q+1, \\
              w(B)C\pi\opi w(\pi, A^{-1}), &r<q+1. 
\end{cases}
\end{split}
\]
\end{lemma}

\begin{proof} The groups for \(A_q^{-1}C_r\) and \(B_q^{-1}C_r\) are
inverses of cases covered in Lemmas \ref{MoveALeft} and
\ref{MoveBLeft}.  The first line for \(C_q^{-1}C_r\) is handled much
as in the proof of the case of \(C_q^{-1}A_r\) in Lemma
\ref{MoveALeft}  and the third line for \(C_q^{-1}C_r\) is the
inverse of the first line.  The second line of \(\pi_qC_r\) is done
by 
\[
\begin{split}
\pi_qC_r
&=
\pi_qC_{q+2}
(B_{q+1}\pi_{q+2}A_{q+1}^{-1})(B_{q}\pi_{q+1}A_{q}^{-1})
\cdots (B_{r}\pi_{r+1}A_{r}^{-1})
\\ &=
C_{q+2}\pi_q
(B_{q+1}\pi_{q+2}A_{q+1}^{-1})(B_{q}\pi_{q+1}A_{q}^{-1})
\cdots (B_{r}\pi_{r+1}A_{r}^{-1}).
\end{split}
\]
The second line for \(\opi_qC_r\) is done by 
\[ \begin{split} \opi_qC_r &=
A_qA_{q+1}\cdots A_{r-2}\opi_{r-1} \pi_{r-2}\pi_{r-3} \cdots \pi_q
C_r \\ &= A_qA_{q+1}\cdots A_{r-2}\opi_{r-1} C_r \pi_{r-2}\pi_{r-3}
\cdots \pi_q \\ &= A_qA_{q+1}\cdots A_{r-2} B_{r-1}\opi_r \pi_{r-1}
\pi_{r-2}\pi_{r-3} \cdots \pi_q. \end{split}\]
The third line for \(\opi_qC_r\) is the worst.  As preparation, we
write 
\[\begin{split} C_r &= C_{q+1}(B_q\pi_{q+1}A_q^{-1})(B_{q-1}\pi_q
A_{q-1}^{-1}) \cdots (B_r\pi_{r+1}A_r^{-1}) \\ &=
C_{q+1}(B_qB_{q-1} \cdots B_r)(\pi_{2q-r+1}A_{2q-r}^{-1}) 
(\pi_{2q-r-1}A_{2q-r-2}^{-1}) \cdots (\pi_{r+1}A_r^{-1}) 
\end{split}\]
which follows from the first line of Lemma \ref{Raises} and from 
the
relations \tref{TwoVRelA}--\tref{TwoVRelB}.  Now we can
write
\[ \begin{split} \opi_q&C_r \\ &= 
\opi_q
C_{q+1}(B_qB_{q-1} \cdots B_r)(\pi_{2q-r+1}A_{2q-r}^{-1})
(\pi_{2q-r-1}A_{2q-r-2}^{-1}) \cdots (\pi_{r+1}A_r^{-1})
\\ &=
B_q\opi_{q+1}\pi_qB_qB_{q-1} \cdots B_r(\pi_{2q-r+1}A_{2q-r}^{-1})
(\pi_{2q-r-1}A_{2q-r-2}^{-1}) \cdots (\pi_{r+1}A_r^{-1}) 
\\ &=
B_q\opi_{q+1}B_{q+1}\pi_q\pi_{q+1}
B_{q-1} \cdots B_r
\\ &\qquad\qquad
(\pi_{2q-r+1}A_{2q-r}^{-1})
(\pi_{2q-r-1}A_{2q-r-2}^{-1}) \cdots (\pi_{r+1}A_r^{-1}) 
\\ &=
B_q\opi_{q+1}B_{q+1}
B_{q-1} \cdots B_r
\pi_{2q-r}\pi_{2q-r+1}
\\ &\qquad\qquad 
(\pi_{2q-r+1}A_{2q-r}^{-1})
(\pi_{2q-r-1}A_{2q-r-2}^{-1}) \cdots (\pi_{r+1}A_r^{-1}) 
\\ &=
B_q\opi_{q+1}B_{q+1}
B_{q-1} \cdots B_r
\pi_{2q-r}A_{2q-r}^{-1}
\\ &\qquad\qquad 
(\pi_{2q-r-1}A_{2q-r-2}^{-1})
(\pi_{2q-r-3}A_{2q-r-4}^{-1}) \cdots (\pi_{r+1}A_r^{-1}) 
\\ &=
B_qC_{q+2}\pi_{q+1}\opi_{q+2}
B_{q-1} \cdots B_r
\pi_{2q-r}A_{2q-r}^{-1}
\\ &\qquad\qquad 
(\pi_{2q-r-1}A_{2q-r-2}^{-1})
(\pi_{2q-r-3}A_{2q-r-4}^{-1}) \cdots (\pi_{r+1}A_r^{-1}) 
\\ &=
B_q
B_{q-1} \cdots B_r
C_{2q-r+2}\pi_{2q-r+1}\opi_{2q-r+2}
\pi_{2q-r}A_{2q-r}^{-1}
\\ &\qquad\qquad 
(\pi_{2q-r-1}A_{2q-r-2}^{-1})
(\pi_{2q-r-3}A_{2q-r-4}^{-1}) \cdots (\pi_{r+1}A_r^{-1}).
\end{split}\]
\end{proof}

\subsection{The \protect\(LMR\protect\) form, Part II} 

We can now add \(C\) and \(C^{-1}\) to the list of generators that
we can handle.

\begin{lemma}\mylabel{LMRABCpi} Let \(w\) be of the form
\(w(A,B,C,\pi,A^{-1}, B^{-1},C^{-1})\).  Then \(w\sim LMR\) where
\(L\) and \(R^{-1}\) are words of the form \(w(A,B,C)\) and \(M\) is
of the form \(w(\pi)\).  Further the number of appearances of \(C\)
in \(L\) will be no larger than the number of appearances of \(C\)
in \(w\) and the number of appearances of \(C^{-1}\) in \(R\) will
be no larger than the number of appearances of \(C^{-1}\) in \(w\).
\end{lemma}

\begin{proof} Let \(w\) be a word of form \(w(A,B,C,\pi,A^{-1},
B^{-1},C^{-1})\).  We will deal in syllables of \(w\).  In this
proof a syllable will be a maximal subword of \(w\) of form
\(w(A,B,\pi,A^{-1}, B^{-1})\).  Thus \(w\) is an alternation of
syllables and words of form \(w(C,C^{-1})\).  We will alter the word
\(w\) using the information in Lemma \ref{MoveCLeft}.  At each
stage, we can assume that each syllable is in the \(LMR\) form of
Lemma \ref{LMRABpi}.

All relations that we will use in this argument will not raise the
number of appearances of \(C\) and \(C^{-1}\) and the last sentence
of the lemma will follow from inspection the arguments.

Assume first that \(w\) has a \(C^{-1}\) that appears somewhere to
the left of an appearance of \(C\) in \(w\).  If these are adjacent,
then the number of appearances of \(C\) and \(C^{-1}\) can be
lowered by using the third group from Lemma \ref{MoveCLeft}.  If
there is no such adjacency, then there is a syllable with \(C^{-1}\)
on the left and \(C\) on the right.  Using the first, second and
fourth groups from Lemma \ref{MoveCLeft}, the \(C\) to the right can
be moved over the maximal subword of form \(w(\pi, A^{-1}, B^{-1})\)
of the syllable at the expense of making the syllable to the right
of the \(C\) more complicated.  Using the inverses of the same
groups from Lemma \ref{MoveCLeft}, we can move the \(C^{-1}\) on the
left over the remaining part of the syllable which now has the form
\(w(A,B)\) at the expense of making the syllable to the left of the
\(C^{-1}\) more complicated.  Now the \(C\) and \(C^{-1}\) are
adjacent and can be eliminated as before.

Thus we can assume that all appearances of \(C\) in \(w\) are to the
left of all appearances of \(C^{-1}\).

Let \(p\) be the largest prefix of \(w\) and let \(s\) be the
largest suffix of \(w\) with \(p\) and \(s^{-1}\) both of form
\(w(A,B,C)\).  We are done if we can get all appearances of \(C\) in
\(p\) and all appearances of \(C^{-1}\) in \(s\).

Consider the leftmost appearances of \(C\) that is not in \(p\).  It
is separated from \(p\) by a syllable.  This syllable must be in
\(LMR\) form as in Lemma \ref{LMRABpi}.  Further this syllable must
also be of form \(w(\pi, A^{-1}, B^{-1})\) since the \(L\) part will
be absorbed by \(p\).  As above, we use the inverses of groups one,
two and four from Lemma \ref{MoveCLeft} to move the \(C\) past all
letters in the syllable at the expense of making the syllable to the
right of the \(C\) more complex.  Inductively we get all appearances
of \(C\) in \(p\).  The appearances of \(C^{-1}\) are handled
similarly.  \end{proof}

\subsection{The \protect\(LMR\protect\) form, Part III} 

We can now add \(\opi\) to the list of generators that
we can handle.

\begin{lemma}\mylabel{LMRABCpiopi} Let \(w\) be a word in
\(\Sigma_s\) of Section \ref{ConventionSec}.  Then \(w\sim LMR\)
where \(L\) and \(R^{-1}\) are words of the form \(w(A,B,C)\) and
\(M\) is of the form \(w(\pi,\opi)\).  \end{lemma}

\begin{proof} We sketch the argument.  We will only use relations
that do not alter the number of appearances of the \(\opi\) in a
word.  We will exploit the fact that the interchange rules of Lemma
\ref{MoveALeft} are the least complex.

We are concerned with syllables that are maximal of the form
\(w(\pi,\opi)\).  If there is more than one such syllable in a word
\(w\), then there are two \(s_1\) and \(s_2\) that are separated by
a word in the \(LMR\) form of Lemma \ref{LMRABCpi} giving a subword
\(s_1LMRs_2\).  In \(L\) we find generators \(A\), \(B\) and \(C\).

Using the fourth and fifth groups from Lemmas \ref{MoveALeft} and
\ref{MoveBLeft}, we can pass appearances of \(A\) and \(B\) from
\(L\) over a single appearance of \(\opi\) in \(s_1\) at the expense
of introducing more complicated expressions to the left of the
\(\opi\) and copies of \(\pi\) to the right of the \(\opi\).

From the fourth and fifth groups from Lemma \ref{MoveCLeft}, we can
move a copy of \(C\) to the left at greater expense.  Appearances of
\(A^{-1}\) and \(B\) will be made to the right of the \(\pi\) or
\(\opi\) that is crossed over.  When put in the \(LMR\) form of
Lemma \ref{LMRABpi} we get copies of \(A^{-1}\) that have to move to
the right and copies of \(B\) that have to move to the left.  Using
Lemma \ref{LMRABCpi} and the inverses of the formulas in Lemma
\ref{MoveALeft}, we see that the copies of \(A^{-1}\) can be
migrated completely to the \(R\) part of the altered \(s_1LMRs_2\)
making what is left of the \(LM\) part more complicated, but without
raising the number of appearances of \(C\) that are left in the
\(L\) part.

Thus the appearances of \(C\) in \(L\) can be passed over each
letter in \(s_1\) as well as the (increasing number) of \(A\) and
\(B\) generators between them.  Eventually, \(s_1LMRs_2\) is reduced
to a word of the form \(w(A,B,C)s_3R's_2\) where \(s_3\) is the
altered form of \(s_1\).

Now we apply the inverses of what we have done to pass \(R\) over
\(s_2\).  This will result in the introduction of copies of \(A\)
which will have to pass over \(s_3\).  Eventually, \(s_1LMRs_2\) is
reduced to a word of the form \(w(A,B,C)s_4w(A^{-1}, B^{-1},
C^{-1})\) where \(s_4\) is the combination of the altered form of
\(s_3\) and \(s_2\).  This reduces the number of syllables by one.

We now assume that our original word \(w\) is of the form \(psq\)
where \(p\) and \(q\) are in the \(LMR\) form of Lemma
\ref{LMRABCpi}.  Thus \(w=LMRsL'M'R'\) with the obvious comments.
As before, we pass all of \(R\) over \(s\) and put the right side
again in \(LMR\) form of Lemma \ref{LMRABCpi}.  Now the new \(L'\)
is passed to the left.  All the while extra instances of \(A\) or
\(A^{-1}\) that have to migrate ``the other way'' are handled as
above.  Eventually, we reach our goal.  \end{proof}

\subsection{Improving \protect\(L\protect\) and
\protect\(R\protect\), Part I} We take the first step in getting the
\(L\) and \(R\) parts of \(LMR\) in a more canonical form.  In the
following, we have to allow \(n=-1\) since \(p\) might be the empty
word.  Similarly, we have to allow \(m=-1\).

\begin{lemma}\mylabel{CsInFront} Let \(w\) be a word in \(\Sigma_s\)
of Section \ref{ConventionSec}.  Then \(w\sim LMR\) as in Lemma
\ref{LMRABCpiopi} and in addition, \(L=pq\) where \(p=C_{i_0}
C_{i_1}\cdots C_{i_n}\) with \(n\ge-1\), with \(i_0< i_1< \cdots <
i_n\) and \(q\) is a word of form \(w(A,B)\), and \(R^{-1}=p'q'\)
where \(p'=C_{j_0} C_{j_1}\cdots C_{j_m}\) with \(m\ge-1\), with
\(j_0< j_1< \cdots < j_m\) and \(q'\) is a word of form \(w(A,B)\).
\end{lemma}

\begin{proof} Let \(L\) be as given by Lemma \ref{LMRABCpiopi}.  Let
\(p\) be the longest (possibly empty) prefix of \(L\) of the form
\(L=C_{i_0} C_{i_1}\cdots C_{i_n}\) with \(i_0< i_1< \cdots < i_n\)
and let \(r\) be the remainder of \(L\) in that \(L=pr\).  Let
\(C_j\) be the leftmost appearance of the generator \(C\) in \(s\).
We have \(L\) in the correct form if there is none.

We write \(s=uC_jv\).  Using Lemma \ref{Raises}, we can raise the
subscript of \(C_j\) as far as we like at the expense of introducing
a word of form \(w(B, \pi, A^{-1})\) before \(v\).  As in the proof
of Lemma \ref{LMRABCpiopi}, we can move the appearances of
\(A^{-1}\) past \(vM\) without raising the number of appearances of
\(C\) in \(v\).

Using the previous paragraph, we raise the subscript of \(C_j\) so
that it is higher than \(i_n\) plus the maximum of all the
subscripts in \(u\) plus the number of letters in \(u\).  Since
\(u\) is a word of form \(w(A,B)\), we can use \tref{TwoVRelH} to
pass the altered \(C_j\) to the left of each letter in \(u\),
lowering the subscript of \(C\) by one with each application of
\tref{TwoVRelH}.  Our elevation of the subscript guarantees that
\tref{TwoVRelH} applies at each step of this passage and that the
ending subscript will be higher than \(i_n\).  We have \(L\) in the
right form by induction.

We now look at \((LMR)^{-1} = R^{-1}M^{-1}L^{-1}\) and apply what we
have done to \(R^{-1}\).  We get the right form for \(R^{-1}\) at
the expense of adding a word of form \(w(A^{-1})\) to the left of
\(L^{-1}\).  This keeps the correct form for \(L\).  \end{proof}

\subsection{Structure from
\protect\(L\protect\)}\mylabel{LeftExtDefSec}

We will extract structure from \(L\) (and \(R^{-1}\)) assumed to be
in the form from Lemma \ref{CsInFront}.  We will do so inductively,
so we will have to describe the structure before we prove it exists.

We will show that as an element in \(\widehat{2V}\), the pair
representing such an \(L\) will be in the form \((t, v_0^k)\) where
\(k\) is the length of \(t\) and \(t\) is of form \(w(v,h)\).
Further, the word \(t\) will correspond to a forest whose only
non-trivial tree \(T\) is the 0-th tree.  Note that \(v_0^k\) also
corresponds to a forest whose only non-trivial tree is the 0-th
tree.  We call the tree \(T\), the \emph{tree corresponding to}
\(t\).

With \(L\), \(t\), \(T\) and \(k\) as in the previous paragraph,
\(L=(tv_0^j, v_0^{j+k})\) also holds for any \(j\ge0\).  From our
methods of building forests from an element of \(\Pi\) described in
Section \ref{BuildForSec} and from the fact given in Lemma
\ref{WordsInVHNumber} that the leaf numbering of a forest
corresponding to a word of form \(w(v,h)\) is the standard
left-right numbering, we know that appending \(v_0^j\) to the right
of \(t\) just adds a caret to the leftmost leaf of the forest for
\(t\) repeatedly \(j\) times.  This is pictured below.

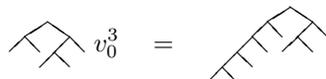
\[
\xy
(0,-1); (2,1)**@{-}; (5,3)**@{-}; (8,1)**@{-}; (4,-3)**@{-};
(2,1); (4,-1)**@{-};   (8,1); (10,-1)**@{-};
(6,-1); (8,-3)**@{-};
\endxy
\,\,v_0^3
\quad
=
\quad
\xy
(0,1); (2,3)**@{-}; (5,5)**@{-}; (8,3)**@{-}; (4,-1)**@{-};
(2,3); (4,1)**@{-};   (8,3); (10,1)**@{-};
(6,1); (8,-1)**@{-};
(2,-1); (0,1)**@{-}; (-6,-5)**@{-};
(-2,-1); (0,-3)**@{-}; (-4,-3); (-2,-5)**@{-};
\endxy
\]

Thus the tree corresponding to \(tv_0^j\) is obtained from \(T\) by
adding a caret with label \(v\) to the leftmost leaf of \(T\)
exactly \(j\) times.  We refer to this as an \emph{extension of
\(T\) to the left}.  Of course, the same thing happens in the
passage from \(v_0^k\) to \(v_0^{k+j}\).

\subsection{The right-left leaf order}

We saw in Section \ref{BuildForSec} how the leaf numbering and the
letter subscripts cooperated in telling where the next caret is to
be attached.  We also saw in Lemma \ref{WordsInVHNumber} that the
leaf numbering of a forest for a word of form \(w(v,h)\) is the
standard left-right numbering of the leaves.

We will discover that the left-right leaf numbering will not
cooperate well with the subscripts of the letters in a word of form
\(w(A,B,C)\).  However, a right-left numbering does.  Since all the
information from a word of form \(w(A,B,C)\) is concentrated in a
single tree, we will only discuss trees here.  Later, we will extend
the discussion to forests.

Given a tree \(T\), we will refer to two numberings of the leaves.
If the tree has \(k\) carets (\(k\) internal vertices), it will have
\(k+1\) leaves which can be numbered from \(0\) through \(k\).  The
\emph{left-right numbering} and \emph{right-left numbering} should
be self descriptive, but we make sure by pointing out that the
following is true of each vertex in \(T\) in the left-right
numbering: all leaves below the left child are numbered less than
all the leaves below the right child.  For the right-left numbering,
the phrase ``less than'' is replaced by ``greater than.''  Note that
for each leaf of the tree, the numbers from the two numberings will
add up to \(k\).

\subsection{Extending the right-left leaf order}

Let \(T\) be a tree corresponding to a word \(t\) of form \(w(v,h)\)
that arrises from an \(L\) from Lemma \ref{CsInFront}.  An extension
\(T'\) of \(T\) to the left corresponds to \(tv_0^j\) for some
\(j\ge0\).  The extra carets that make \(T'\) from \(T\) are
constantly added to the leftmost leaf.  Thus the right-left leaf
order on \(T\) carries over to those leaves that \(T\) and \(T'\)
have in common.  On \(T\), this consists of all leaves of \(T\)
except the leftmost.  This observation will be used repeatedly in
what follows.

\subsection{Building a tree from \protect\(L\protect\)}

From Lemma \ref{CsInFront}, we are motivated to study words such as
\(C_{i_0}C_{i_1}\cdots C_{i_n}w(A,B)\) with \(i_0<i_1<\cdots<i_n\).
Later we will be obliged to apply relations \tref{TwoVRelI} and
\tref{TwoVRelK} to such words which will bring in elements of
\(\{\pi_i \mid i\in\N\}\).  With a little work and Lemma
\ref{LMRABpi}, we will be able to move appearances of \(\{\pi_i \mid
i\in\N\}\) to the end of the words.  This briefly justifies our
concentration on the words that appear in the next few lemmas.

We start without any appearances of \(\pi_i\).  Let
\(L=C_{i_0}C_{i_1}\cdots C_{i_n}w(A,B)\) with
\(i_0<i_1<\cdots<i_n\).  Let \(l\) be the length of \(L\).  For each
\(j\) with \(0\le j\le l\), let \(p_j\) be the prefix of \(L\) of
length \(j\).  We will show that \(L\) corresponds to an element of
the form \((t, v_0^k)\) in \(\widehat{2V}\) and we want to describe
\(t\) and the tree \(T\) corresponding to \(t\).  We will do so
inductively by describing these items for each \(p_j\) and how they
are obtained from the corresponding items for \(p_{j-1}\).  The
prefix \(p_0\) is the empty string and its element of
\(\widehat{2V}\) is \((v_0^0, v_0^0)\) and its tree is the trivial
tree.  We write \(p_j=(t_j, v_0^{k_j})\) and the tree for \(t_j\) is
\(T_j\).  Note that for \(1\le j\le n+1\), we have
\(p_j=C_{i_0}\cdots C_{i_{j-1}}\).  We now give the inductive
lemmas.

\begin{lemma}\mylabel{AddAC} If we take the notation and assumptions
of the previous paragraph and restrict \(j\) so that \(1\le j\le
n+1\), then we have {\claimenum \item \(k_j=i_{j-1}+1\), \item
\(t_j=t_{j-1}v_0^dh_0\) where \(d=i_{j-1}-i_{j-2}-1\), and \item
\(T_j\) is obtained from \(T_{j-1}\) by attaching a caret labeled
\(h\) to the leaf numbered \(i_{j-1}\) in the right-left leaf order
in the smallest left extension of \(T_{j-1}\) that has a leaf
numbered \(i_{j-1}\) in the right-left leaf order. \claimenumend}
\end{lemma}

\begin{proof} Item (a) follows from (b) by induction and item (c)
follows from (b) directly.  Thus we must show (b).

Let \(m=i_{j-1}\) and \(n=k_{j-1}\) for typographical reasons.  We
have \(m>i_{j-2}\) by assumption and \(i_{j-2}=n-1\) by induction,
so \(m\ge n\).  We set \[d=m-n = i_{j-1}-i_{j-2}-1.\] From Section
\ref{TwoVGenSec}, we have \(C_m=(v_0^mh_0,v_0^{m+1})\).  Now
\[\begin{split} p_j= p_{j-1}C_{i_{j-1}} = p_{j-1}C_m &= (t_{j-1},
v_0^n)(v_0^mh_0, v_0^{m+1}) \\ &= (t_{j-1}v_0^dh_0,
v_0^nv_0^dh_0)(v_0^mh_0, v_0^{m+1}) \\ &= (t_{j-1}v_0^dh_0,
v_0^mh_0)(v_0^mh_0, v_0^{m+1}) \\ &=(t_{j-1}v_0^dh_0, v_0^{m+1})
\end{split}\] which is what we needed to show. \end{proof}

\begin{lemma}\mylabel{AddAnX} Take the notation and assumptions of
the paragraph before Lemma \ref{AddAC} and restrict \(j\) so that
\(n+1< j\le l\).  Let \(X_i\) be such that \(p_j=p_{j-1}X_i\) with
\(X\) one of \(\{A,B\}\).  Let \(x=v\) if \(X=A\) and \(x=h\) if
\(X=B\).  Let \(n=k_{j-1}\).  Then \item \[p_{j-1}X_i =
\begin{cases} (t_{j-1}v_0^{i+1-n}x_1, v_0^{i+1}), &n\le i+1, \\
(t_{j-1}x_{n-i}, v_0^{n+1}), &n>i+1. \end{cases}\] and \(T_j\) is
obtained from \(T_{j-1}\) by attaching a caret labeled \(x\) to the
leaf numbered \(i\) in the right-left leaf order in the smallest
left extension of \(T_{j-1}\) that has at least \(i+2\) leaves.
\end{lemma}

\begin{proof}  From Section
\ref{TwoVGenSec}, we have \(X_i=(v_0^{i+1}x_1,v_0^{i+2})\).
If \(n\le i+1\), then we have \[\begin{split}
p_{j-1}X_i &= (t_{j-1},v_0^n)(v_0^{i+1}x_1, v_0^{i+2}) \\ &=
(t_{j-1}v_0^{i+1-n}x_1, v_0^{i+1}x_1) (v_0^{i+1}x_1, v_0^{i+2}) \\
&= (t_{j-1}v_0^{i+1-n}x_1, v_0^{i+2}). \end{split}\] If \(n>i+1\),
then we have \[\begin{split} p_{j-1}X_i &=
(t_{j-1},v_0^n)(v_0^{i+1}x_1, v_0^{i+2}) \\ &=(t_{j-1},v_0^n)
(v_0^{i+1}x_1v_0^{n-i-1}, v_0^{n+1}) \\ &= (t_{j-1},v_0^n)
(v_0^nx_{n-i}, v_0^{n+1}) \\ &=(t_{j-1}x_{n-i}, v_0^nx_{n-i})
(v_0^nx_{n-i}, v_0^{n+1}) \\ &= (t_{j-1}x_{n-i},
v_0^{n+1}). \end{split}\] 

When \(n\le i+1\), we are adding a caret with label \(x\) at leaf 1
in the left-right order to a tree with \(i+1\) carets and thus
\(i+2\) leaves numbered from 0 through \(i+1\).  Thus the addition
is at leaf \(i\) in the right-left order and the tree is the
smallest left extension of \(T_{j-1}\) that has at least \(i+2\)
leaves.  When \(n>i+1\), we are adding a caret with label \(x\)
directly to \(T_{j-1}\) at the leaf numbered \(k_{j-1}-i\) in the
left-right order, or the leaf numbered \(i\) in the right-left
order.  Note that in this case, the tree \(T_{j-1}\) already has at
least \(i+2\) leaves.  \end{proof}

We now add appearances of the \(\pi_i\), but don't worry about the
tree structure.

\begin{lemma}\mylabel{AddAPi} Let \(p=C_{i_0}C_{i_1}\cdots
C_{i_n}w(A,B)w(\pi)\) with \(i_0<i_1<\cdots <i_n\) and assume that
\(p=(t,v_0^n)\) where \(t\) is a word of form \(w(v,h,\sigma)\) and
\(n\) is the number of appearances of \(v\) and \(h\) in \(t\).
Then \(p\pi_i = (tv_0^j\sigma_{(n+j-1)-i}, v_0^{n+j})\) where \(j\)
is the smallest value in \(\N\) so that \((n+j-1)-i>0\).  In
particular, \(j=0\) if \((n-1)-i>0\) (equivalently, \(n\ge i+2\)),
and \((n+j-1)-i=1\) if \(n<i+2\).  \end{lemma}

\begin{proof} A calculation similar to that in Lemma \ref{AddAnX}
using \(\pi_i=(v_0^{i+2}\sigma_1, v_0^{i+2})\) gives \[p\pi_i
= \begin{cases} (tv_0^{i+2-n}\sigma_1, v_0^{i+2}), &n< i+2, \\
(t\sigma_{(n-1)-i}, v_0^n), &n\ge i+2. \end{cases}\] The rest is
straightforward.  \end{proof}

\begin{lemma}\mylabel{PisWereAdded} Let \(L=C_{i_0}C_{i_1}\cdots
C_{i_n}w(A,B)\) with \(i_0<i_1<\cdots <i_n\), and let \(p=Lv\) with
\(v\) a word in \(\{\pi_i\mid 0\le i\le k\}\).  Let \(r\in\N\) be
such that \(n+r\ge k+2\).  Then \(L=(t, v_0^n)\) where \(t\) is a
word of form \(w(v,h)\) and \(n\) is the length of \(t\), and
\(p=(tv_0^rs, v_0^{n+r})\) where \(s\) is a word of form
\(w(\sigma)\).  \end{lemma}

\begin{proof} The claim about \(L\) follows from Lemmas \ref{AddAC}
and \ref{AddAnX}.  The claim about \(p\) follows from Lemma
\ref{AddAPi}.  This is seen by noting that if we set \(L=(tv_0^r,
v_0^{n+r})\), then the only case that arises in applying Lemma
\ref{AddAPi} to each letter in \(v\) is the case in which the
\(j\) of that lemma is equal to 0.  \end{proof}

\subsection{The primary tree from
\protect\(L\protect\)}\mylabel{PrimTreeSec}

We show how much flexibility there is in representations of the form
\((t, v_0^k)\) in \(\widehat{2V}\).

\begin{lemma}\mylabel{TreeTrunkReps} Let \((t, v_0^k)=(s, v_0^j)\)
in \(\widehat{2V}\) where \(k\le j\).  Then \(s=tv_0^{j-k}\).
\end{lemma}

\begin{proof} From the structure of \(\widehat{2V}\) as a group of
right fractions of \(\Pi\), there are \(p\) and \(q\) in \(\Pi\) so
that \((tp,v_0^k p)=(sq,v_0^j q)\) as pairs, giving \(p=v_0^{j-k}q\)
from the cancellativity of \(\Pi\).  The claim follows from
\(tv_0^{j-k}q=tp=sq\).  \end{proof}

This immediately gives the following.

\begin{lemma}\mylabel{LAsWordInVH} Let \(L=C_{i_0}C_{i_1}\cdots
C_{i_n}w(A,B)\) with \(i_0<i_1<\cdots<i_n\).  Let \(k\) be the
smallest in \(\N\) so that \(L=(t, v_0^k)\) as an element in
\(\widehat{2V}\) where \(t\) is of form \(w(v,h)\) and \(k\) is the
length of \(t\).  Let \(L=(s, v_0^j)\) where \(s\) is of form
\(w(v,h)\) and \(j\) is the length of \(s\).  Then \(s=tv_0^{j-k}\).
\end{lemma}

With \(L\), \(k\) and \(t\) as in Lemma \ref{LAsWordInVH}, we call
the tree \(T\) corresponding to \(t\), the \emph{primary tree for
the word \(L\)}.  If \(L=(s, v_0^j)\) is any other representation of
\(L\) with \(s\) a word of form \(w(v,h)\) and \(j\) equal to the
length of \(s\), then we know that the tree corresponding to \(s\)
is an extension of \(T\) to the left.  Since
\(L=(tv_0^n,v_0^{k+n})\) is a valid representation of \(L\) of the
correct form for any \(n\in\N\), we know that all extensions of \(T\)
to the left can show up in this way.  We call such an extension with
a \emph{secondary tree for \(L\)} even when \(n=0\).  Thus the
primary tree for \(L\) is also a secondary tree for \(L\).

\subsection{The stucture of the primary tree for
\protect\(L\protect\)} \mylabel{LTreeStructSec}

We consider a word \[L=C_{i_0}C_{i_1}\cdots C_{i_n} X_{i_{n+1}}
\cdots X_{i_{l-1}}\] where \(i_0<i_1<\cdots i_n\) and where each
\(X\) comes separately from \(\{A,B\}\).  The primary tree for \(L\)
will be described as a particularly simple tree with a finite forest
attached.  The right-left leaf order will be used throughout.  We
must state how this order extends to forests, since we use one of
two obvious choices and have to be explicit as to which.

The attached forest is finite in that it has finitely many trees.
The right-left leaf order of a finite forest \(F\) with trees
\(F_i\) with \(0\le i\le q\) has its leaves numbered consecutively
starting from \(0\) with all leaves in \(F_i\) numbered above those
in \(F_j\) whenever \(i>j\) and the numbering of the leaves in any
one \(F_i\) following the right-left leaf order.  This is best
pictured if the forest is drawn with \(F_0\) the rightmost tree and
\(F_q\) the leftmost.  This is the reverse of the usual picture.

The particularly simple tree that our forest is attached to will
correspond to the subword \(C_{i_0}C_{i_1}\cdots C_{i_n}\) and will
have the form pictured below.
\mymargin{Trunk}
\begin{equation}\label{Trunk}
\begin{split}
\xy
(-6,-6); (0,0)**@{-}; (2,-2)**@{-};
(-4,-4); (-2,-6)**@{-}; (-2,-2); (0,-4)**@{-};
(-18,-18); (-12,-12)**@{-}; (-10,-14)**@{-};
(-16,-16); (-14,-18)**@{-}; (-14,-14); (-12,-16)**@{-};
(-10,-10)*{\cdot}; (-9,-9)*{\cdot}; (-8,-8)*{\cdot};
\endxy
\end{split}
\end{equation}
The tree in \tref{Trunk} is the result of any word of the form
\(a_0a_1\cdots a_n\) where each \(a_i\) is from \(\{v_0, h_0\}\).
Note that the subscript \(i\) of  \(a_i\) is not part of the symbol
that \(a_i\) represents.  The tree in \tref{Trunk} will be referred
to as a \emph{trunk}.

If \(\Lambda \) is a trunk with \(m\) carets and \(m+1\) leaves with
the leaves numbered from \(0\) through \(m\) in the right-left
order, and \(F\) is a finite forest with \(m\) trees, then we
combine \(\Lambda \) with \(F\) to produce a tree \(T\) by attaching
each \(F_i\) with \(0\le i< m\) to the leaf in \(\Lambda \) numbered
\(i\) in the right-left order.  It is deliberate that we never use
the leaf in \(\Lambda \) that is numbered \(m\) (the leftmost leaf).

The forest \(F\) will be built from \( X_{i_{n+1}} \cdots
X_{i_{l-1}}\) much as forests are built from words in Section
\ref{BuildForSec} with a few differences.  As in Section
\ref{BuildForSec}, we build the forest caret by caret as we build
the word from left to right letter by letter.  The forests we build
here will end up with the right-left ordering on the leaves.  We
start with the trivial forest and then for each \(A_i\), we add a
caret to leaf \(i\) with label \(v\), we keep the numbering on all
leaves with number less than \(i\), we increase by 1 the numbers of
all leaves with number greater than \(i\), we number the new left
leaf \(i+1\), and we number the new right leaf \(i\).  The numbering
of the new leaves is different from the scheme in Section
\ref{BuildForSec}.  For each \(B_i\), we do exactly the same thing,
except the label of the new caret is \(h\).

\begin{lemma}\mylabel{LTreeStruct} Let \(L=C_{i_0}C_{i_1}\cdots
C_{i_n} X_{i_{n+1}} \cdots X_{i_{l-1}}\) where \(i_0<i_1<\cdots
i_n\) and where each \(X\) comes separately from \(\{A,B\}\).  Let
\(m\) equal the maximum of \[\{i_j+n+2-j\mid n+1\le j\le
l-1\}\cup\{i_n+1\}.\] Then \(L\) can be represented as \(L=(t,
v_0^k)\) where \(t\) is of form \(w(v,h)\) and \(k\) is the length
of \(t\), so that \(k=m+l-n\), and so that the tree \(T\) for \(t\)
is the primary tree for \(L\) and is described as follows.  The tree
\(T\) consists of a trunk \(\Lambda \) with a finite forest \(F\)
attached.  The trunk \(\Lambda \) has \(m\) carets and \(m+1\)
leaves numbered \(0\) through \(m\) in the right-left order.  If the
carets in \(\Lambda \) are numbered from 0 starting at the top, then
the label of the \(i\)-th caret is \(h\) if \(i\) is in \(\{i_0,i_1,
\ldots\, i_n\}\) and \(v\) otherwise.  The forest \(F\) is built
from the word \( X_{i_{n+1}} \cdots X_{i_{l-1}}\) as described just
prior to this statement.  \end{lemma}

\begin{proof} We can discuss the tree \(T\) as built from \(L\)
letter by letter because of Lemmas \ref{AddAC} and \ref{AddAnX}.
The subtree \(T'\) coming from the prefix \(C_{i_0}C_{i_1}\cdots
C_{i_n}\) that we get from Lemma \ref{AddAC} is as \(\Lambda \) is
described (including the labeling) in the statement above, except
that the number of carets in \(T'\) will only by \(i_n+1\).  The
trunk \(\Lambda \) is an extension of \(T'\) to the left so as to
have \(m\) carets and highest leaf number \(m\).

To add \(X_i\) to an existing tree, we need to create a left
extension of the tree if the tree has no leaf numbered \(i+1\).  Let
us assume that we can build the tree from \(\Lambda \) as described
in the statement with no extra left extensions needed through
\(X_{i_{j-1}}\).  The tree for \(C_{i_0}C_{i_1}\cdots C_{i_n}
X_{i_{n+1}} \cdots X_{i_{j-1}}\) will have highest leaf number equal
to \(m+j-(n+1)\).  The next letter to be treated will be \(X_{i_j}\)
and our hypothesis dictates that \(m\ge i_j+n+2-j\) or
\(m+j-(n+1)\ge i_j+1\).  Thus the caret for \(X_{i_j}\) can be added
without further extension.  The forest \(F\) is simply the tree
\(T\) for \(L\) as built from Lemmas \ref{AddAC} and \ref{AddAnX}
with the trunk \(\Lambda \) removed where the trunk \(\Lambda \)
consists of all carets reachable from the root by repeatedly going
to the left child.

The number of carets in \(T\) is \(m\) plus the number of carets in
\(F\), and the number of carets in \(F\) is \(l-n\).  The number of
carets in \(T\) must be the length of the word \(t\), so
\(k=m+l-n\).

Note that either \(m+j-(n+1)=i_j+1\) for some \(j\) or \(m=i_n+1\).
In the first case, the bottom caret in the trunk \(\Lambda \) has a
caret attached to one of its leaves.  In the second case, the bottom
caret of \(\Lambda \) has label \(h\).  If the tree \(T\) is not the
primary tree for \(L\), then it is an extension to the left of
another tree by Lemma \ref{LAsWordInVH}.  This is impossible by the
remarks we have just made.  \end{proof}

Note that \(m\) must come out to be at least \(1\) in Lemma
\ref{LTreeStruct}.  This applies even if \(n=-1\) which happens if
there is no appearance of any \(C_i\) in \(L\).  Thus the trunk
\(\Lambda \) is never empty.

\begin{lemma}\mylabel{LTreeStructCor} Let \(L=C_{i_0}C_{i_1}\cdots
C_{i_n}w(A,B)\) with \(i_0<i_1<\cdots<i_n\) and let \(L=(s, v_0^k)\)
where \(s\) is a word of form \(w(v,h)\) and \(k\) is the length of
\(s\).  Then the tree \(T'\) for \(s\) is a secondary tree for
\(L\).  Further, \(T'\) is an extension to the left of the primary
tree \(T\) for \(L\) and has the same description as \(T\) as given
in Lemma \ref{LTreeStruct} except that the trunk for \(T'\) is an
extension to the left of the trunk for \(\Lambda \) in Lemma
\ref{LTreeStruct}.  \end{lemma}

\begin{proof}  This follows immediately from Lemmas
\ref{LAsWordInVH} and \ref{LTreeStruct}, the definitions in
Section \ref{PrimTreeSec}, and the definition of extension to the
left as found in Section \ref{LeftExtDefSec}.  \end{proof}

\subsection{Improving \protect\(M\protect\)} We now take care of the
\(M\) part of the \(LMR\) form.

\begin{lemma}\mylabel{ImprovM} Let \(w\) be a word in \(\Sigma_s\)
of Section \ref{ConventionSec}.  Then \(w\sim LMR\) as in Lemma
\ref{CsInFront} with \(L=(s, v_0^m)\) and \(R^{-1}=(t, v_0^n)\) so
that \(s\) and \(t\) are words of form \(w(v,h)\), \(m\) and \(n\)
are, respectively, the lengths of \(s\) and \(t\), and so that there
is a \(p\ge\max\{m,n\}\) so that \(M\) is a word in \(\{\opi_{p-1},
\pi_i\mid i\le p-2\}\).  \end{lemma}

\begin{proof} This is a direct consequence of Lemmas 4.7 and 4.11
of \cite{brin:bv3} and the definitions made in \cite{brin:bv3} just
before those lemmas.  The cited lemmas of \cite{brin:bv3} apply since
they are about a group called \(BV\) in \cite{brin:bv3}
presented by a generating set \(\{v_i, \pi_i,
\opi_i\mid i \in\N\}\) and a set of relations
given in Lemma 4.2 of
\cite{brin:bv3}
that are all seen to hold in our setting when the
generators of \cite{brin:bv3} are mapped to generators of
\(\Sigma_s\) under the mapping \(v_i\mapsto A_i\), \(\pi_i\mapsto
\pi_i\) and \(\opi_i\mapsto \opi_i\).  The relations \(\pi_i^2=1\)
and \(\opi_i^2=1\) of \(V\) are not relations of \(BV\) 
and as a consequence the relations of 
\(BV\) in \cite{brin:bv3} mention both positive and
negative powers of the \(\pi_i\) and \(\opi_i\).  However, these all
reduce to relations in \tref{TwoVRelA}--\tref{TwoVRelQ} because
\(\pi_i^2=1\) and \(\opi_i^2=1\) are assumed here.  The improvements
of Lemmas 4.7 and 14.11 of \cite{brin:bv3} are gained at the expense
of introducing (after the translation \(v\mapsto A\)) a word of form
\(w(A)\) to the left of \(M\) and a word of form \(w(A^{-1})\) to
the right of \(M\).  The control of the powers of \(v_0\)
(corresponding to \(\lambda_0\) in \cite{brin:bv3}) carries over to
our setting and does not disturb the forms of \(L\) and \(R^{-1}\)
guaranteed by Lemma \ref{CsInFront}.  \end{proof}

In the next two lemmas, it will be more convenient to express
elements of \(\widehat{2V}\) as \(PQ^{-1}\) with \(P\) and \(Q\)
from \(\Pi\) rather than \((P,Q)\).  

The next lemma will be used in analyzing not only \(M\), but also
\(L\) and \(R\).  It is lifted from the proof of Proposition 4.13
of \cite{brin:bv3}.  It gives properties about a certain translation
function and its inverse.  We need some definitions.  Let \(M\) be a
word in \(\{\opi_{p-1}, \pi_i\mid i\le p-2\}\) written as
\mymargin{WordInPiOpi}\begin{equation}\label{WordInPiOpi} M =
Y_{i_1}Y_{i_2}\cdots Y_{i_q} \end{equation} where \(0\le i_j\le p-1
\) and \mymargin{SubsToPis}\begin{equation}\label{SubsToPis} Y_{i_j}
= \begin{cases} \pi_{i_j}, &i_j<p-1, \\ \opi_{p-1},
&i_j=p-1. \end{cases} \end{equation} We let
\mymargin{WordInInvSigma}\begin{equation}\label{WordInInvSigma}
\Psi_{p}(M)=\sigma_{k_1}\sigma_{k_2}\cdots \sigma_{k_q}
\end{equation} where \(k_j=(p-1)-i_j\) for \(1\le j\le q\). In the
other direction, if \(u\) is a word in \(\{\sigma_i\mid i\le p-1\}\)
written as in the right side of \tref{WordInInvSigma}, then
\(\Psi_{p}^{-1}(u)\) is taken to be \(M\) as in \tref{WordInPiOpi}
with letters interpreted by \tref{SubsToPis} with each
\(i_j=(p-1)-k_j\) for \(1\le j\le q\).  Note that \(\Psi_{p}\) and
\(\Psi_{p}^{-1}\) are truly inverse to each other as transformations
on words.

\begin{lemma}\mylabel{PiToSigma} Let \(M\) be a word in
\(\{\opi_{p-1}, \pi_i\mid i\le p-2\}\) and let \(u=\Psi_{p}(M)\).
Then \(M=v_0^puv_0^{-p}\).  Further if \(u\) can be taken to a word
\(u'\) in \(\{\sigma_j\mid 0\le j\le p-1\}\) by relations
\tref{PiRelB}--\tref{PiRelD} from Lemma \ref{PiRels}, then \(M\) can
be taken to a word \(M'\) in \(\{\opi_{p-1}, \pi_i\mid i\le p-2\}\)
by relations \tref{TwoVRelL}--\tref{TwoVRelQ} so that
\(u'=\Psi_{p}(M')\).  \end{lemma}

\begin{proof} That \(M=v_0^puv_0^{-p}\) follows from 
\[\opi_{p-1} = (v_0^p\sigma_0, v_0^p)=v_0^p\sigma_0v_0^{-p}
= v_0^p\sigma_{(p-1)-(p-1)} v_0^{-p}
\]
 and
\[\pi_i =(v_0^{i+2}\sigma_1, v_0^{i+2}) =
(v_0^{i+2}\sigma_1 v_0^{p-i-2}, v_0^p) = (v_0^p\sigma_{(p-1)-i},
v_0^p) =v_0^p\sigma_{(p-1)-i}v_0^{-p}. \]


The last sentence of the
lemma follows by noting that an application of
\tref{TwoVRelL}--\tref{TwoVRelQ} to \(M\) results in an application
of \tref{PiRelB}--\tref{PiRelD} according to the
following assoctiation 
\[
\begin{split} 
\mathrm{\tref{TwoVRelP}\,\, or\,\, \tref{TwoVRelQ}} 
    &\quad\leftrightarrow\quad
\mathrm{\tref{PiRelB}}, \\
\mathrm{\tref{TwoVRelL}\,\, or\,\, \tref{TwoVRelN}} 
    &\quad\leftrightarrow\quad
\mathrm{\tref{PiRelC}}, \\
\mathrm{\tref{TwoVRelM}\,\, or\,\, \tref{TwoVRelO}} 
    &\quad\leftrightarrow\quad
\mathrm{\tref{PiRelD}},
\end{split}
\]
and conversely.  \end{proof}

\begin{lemma}\mylabel{ImprovMCor} Let \(w\) be a word in
\(\Sigma_s\) of Section \ref{ConventionSec}.  Then \(w\sim LMR\) as
in Lemma \ref{CsInFront}, and further when \(L\), \(M\) and
\(R^{-1}\) are expressed as elements of \(\widehat{2V}\), they are
expressible as \(L=s v_0^{-p}\), \(R^{-1}=tv_0^{-p}\) and
\(M=v_0^puv_0^{-p}\) where \(s\) and \(t\) are words of form
\(w(v,h)\), \(u\) is a word in \(\{\sigma_j\mid 0\le j\le p-1\}\),
and the lengths of \(s\) and \(t\) are both \(p\).  Further, if
\(u\) can be reduced to the trivial word using relations
\tref{PiRelB}--\tref{PiRelD} from Lemma \ref{PiRels}, then \(M\) can
be reduced to the trivial word using relations
\tref{TwoVRelL}--\tref{TwoVRelQ}.  \end{lemma}

\begin{proof} If we take \(w\sim LMR\) as given by Lemma
\ref{ImprovM}, then we let \(L\) be represented by \((sv_0^{p-m},
v_0^p)\) and \(R^{-1}\) be represented by \((tv_0^{p-n}, v_0^p)\).
The rest follows from Lemma \ref{PiToSigma}.
\end{proof}

\subsection{Normalizing the trees from \protect\(L\protect\) and
\protect\(R^{-1}\protect\)}\mylabel{NormLSec}

As pointed out in the proof of Lemma \ref{LMRABpi}, the assignments
\(v_i\mapsto A_i\), \(h_i\mapsto B_i\), \(\sigma_i\mapsto \pi_i\)
extend to a group homomorphism from \(\widehat{2V}\) to \(G\) and a
monoid homomorphism from \(\Pi\) to \(G\).  The forest \(F\) built
in Section \ref{LTreeStructSec} from a word of form \(w(A,B)\) is
the mirror image of the forest that would have been built from the
corresponding word of form \(w(v,h)\) built from \(w(A,B)\) by
replacing each \(A\) by \(v\) and each \(B\) by \(h\).  We will use
these facts to improve the word \(L\) even more.

Let \(L=C_{i_0}C_{i_1}\cdots C_{i_n} X_{i_{n+1}} \cdots
X_{i_{l-1}}\) where \(i_0<i_1<\cdots i_n\) and where each \(X\)
comes separately from \(\{A,B\}\).  Let \(T\) be a
secondary tree for the word \(L\) as defined in Section
\ref{PrimTreeSec}.  The labeling on the tree lets us build a
numbered pattern in \(S_0\) from \(T\) and we can use this pattern
to add secondary labels to \(T\) as in Section \ref{SecondLabelSec}.
We can then define when \(T\) is normalized exactly as is done in
that section.

The purpose of this section is to prove the following.

\begin{prop}\mylabel{NormalLTree} Given the notation, hypothese and
conclusion as expressed in Lemma \ref{ImprovMCor}, we can further
assume that the trees for \(s\) and \(t\) are normalized.
\end{prop}

We will build a reductive proof of Proposition \ref{NormalLTree}
from two lemmas.  The lemmas will be applied alternately to reduce
the pattern of non-normalized vertices.  Thus the hypotheses of each
lemma will be designed to match the conclusions of the other.  This
partly explains their rather strange statements.  In the next lemma,
the double appearance of \(n\) is deliberate.

\begin{lemma}\mylabel{ControlPerms} Let \(L=C_{i_0}C_{i_1}\cdots
C_{i_n}u\) and \(L'=C_{k_0}C_{k_1}\cdots C_{k_n}u'\) where
\(i_0<i_1<\cdots<i_n\), where \(k_0<k_1<\cdots <k_n\), where \(u\)
is a word of form \(w(A,B)\), and where \(u'\) is a word of form
\(w(A,B,\pi)\).  Assume that \(L\) is expressible as \((t,v_0^p)\)
as an element of \(\widehat{2V}\) with \(t\) a word of form
\(w(v,h)\) and \(p\) is the length of \(t\).  Let \(m\) be the
number of carets of the trunk of \(T\) and assume that \(m\ge
k_n+1\).

If \(L\sim L'\), then there is a word \(u''\) of form \(w(A,B)\),
and there is a word \(z\) in \(\{\pi_i\mid i\le p-2\}\) so that
setting \(L_1=C_{k_0}C_{k_1}\cdots C_{k_n}u''\) and \(L_2=L_1z'\)
gives that \(L\sim L_2\) and \(L_1\) is expressible as \((t',
v_0^p)\) with \(y\) a word of form \(w(A,B)\) of length \(p\) so
that the tree \(T'\) for \(t'\) is normalized except possibly at
interior vertices in the trunk of the tree, and so that the trunk of
\(T'\) has \(m\) carets.  \end{lemma}

\begin{proof} The homomorphism from \(\widehat{2V}\) to \(G\)
defined by \(v_i\mapsto A_i\), \(h_i\mapsto B_i\), \(\sigma_i\mapsto
\pi_i\) allows us to write \(u'\sim u''z'\) where \(u''\) is a word
of form \(w(A,B)\), where \(z'\) is a word in \(\{\pi_i\mid
i\in\N\}\), and where the forest \(F\) for \(u''\) as built in
Section \ref{LTreeStructSec} is normalized.  Since the primary tree
for \(L_1=C_{k_0}C_{k_1}\cdots C_{k_n}u''\) is a trunk with \(F\)
attached, we have satisfied the normalization requirements.  The
rest of the argument is devoted to improving \(z'\) and
understanding the structure of trees associated to \(L_1\).

From Lemma \ref{PisWereAdded}, we know that
\(L_1=(\widehat{t},v_0^q)\) where \(\widehat{t}\) is a word of form
\(w(v,h)\) and \(q\) is the length of \(\widehat{t}\), and that
\(L_2=(\widehat{t}v_0^rs, v_0^{q+r})\) where \(r\ge0\) and where
\(s\) is a word of form \(w(\sigma)\).  It is seen from Lemma
\ref{AddAPi} that \(s\) is equal to \(\Psi_{q+r}(z')\).  Since Lemma
\ref{PisWereAdded} allows \(r\) to be any sufficiently large value,
we cover all cases by saying that there is a \(k\in\N\) so that
\(p+k=q+r\) and \[L=(tv_0^k,v_0^{p+k}) = (\widehat{t}v_0^r
s,v_0^{q+r})=L_2\] as elements of \(\widehat{2V}\).  Thus the
numbered patterns in the unit square represented by \(tv_0^k\) and
\(\widehat{t}v_0^rs\) are identical.

Since \(s\) is a word of form \(w(\sigma)\), the only difference
between the numbered pattern for \(\widehat{t}v_0^r\) and
\(\widehat{t}v_0^rs\) is in the numbering.  The unnumbered patterns
for \(\widehat{t}v_0^r\) and \(\widehat{t}v_0^rs\) are identical.
Thus the unnumbered patterns for \(tv_0^k\) and \(\widehat{t}v_0^r\)
are identical.

We consider the vertices of a tree that are reachable from the root
by always going to the left.  We call these the ``left edge
vertices.'' In the tree for \(tv_0^k\) let the left edge vertices be
\(a_0\), \(a_1\), \dots, \(a_b\) reading from the top.  This makes
\(a_0\) the root of the tree and \(a_b\) the only non-interior
vertex among them.  Thus all \(a_i\) with \(0\le i<b\) have labels.
Since the trunk for \(T\) the tree for \(t\) has \(m\) carets, we
know that \(b=m+k\).  We also know that for \(m\le i<b\) the label
for \(a_i\) is \(v\) and that the right child of \(a_i\) is a leaf.

In the tree for \(\widehat{t}v_0^r\), let the left edge vertices be
\(a'_0\), \(a'_1\), \dots, \(a'_c\) reading from the top.  Since
\(m\ge k_n+1\), we know from Lemma \ref{LTreeStruct} that the label
for \(a'_i\) is \(v\) for \(m\le i<c\).

The following fact is elementary from the description in Section
\ref{NumPatFromForSec} of how vertices in a labeled tree correspond
to rectangles in the unit square: \begin{enumerate}\item[\((*)\)]
The rectangle corresponding to a left edge vertex depends only on
two numbers; the number of left edge vertices above it and the
number of left edge vertices above it with label
\(h\).\end{enumerate}

From \((*)\), from Lemma \ref{LTreeStruct}, and from the fact that
there are \(n+1\) appearances of a \(C\) in both \(L\) and \(L'\),
we know that the rectangle \(R\) corresponding to \(a_m\) is
identical to the rectangle corresponding to \(a'_m\).  Since \(R\)
is divided \(k\) times vertically according to the carets below
\(a_m\) in the tree for \(tv_0^k\), it must be divided in exactly
the same way by the tree for \(\widehat{t}v_0^r\).  Thus the tree
below \(a'_m\) in the tree for \(\widehat{t}v_0^r\) must consist of
an extension to the left by \(k\) carets all labeled \(v\).  From
this we know that \(r\ge k\).

If we give the left-right leaf numbering to the trees for \(tv_0^k\)
and \(\widehat{t}v_0^r\) (this is the leaf numbering that works with
the subscripts of the letters in \(\{v,h\}\)), then the carets below
\(a_m\) and \(a'_m\) appear as follows in both of these trees.

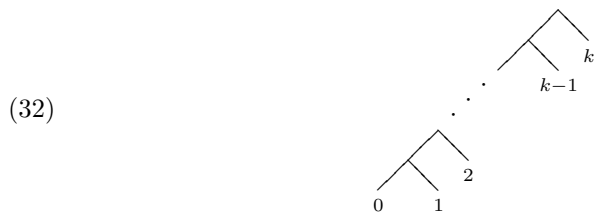
\mymargin{EndOfTrunk}\begin{equation}\label{EndOfTrunk}
\begin{split}
\xy
(0,0); (8,8)**@{-}; (12,4)**@{-}; (4,4); (8,0)**@{-};
(10,10)*{\cdot}; (12,12)*{\cdot}; (14,14)*{\cdot};
(16,16); (24,24)**@{-}; (28,20)**@{-}; (20,20); (24,16)**@{-};
(0,-2)*{\scriptstyle{0}};
(8,-2)*{\scriptstyle{1}};
(12,2)*{\scriptstyle{2}};
(28,18)*{\scriptstyle{k}};
(24,14)*{\scriptstyle{k-1}};
\endxy
\end{split}
\end{equation}

The last \(k\) carets in \(v_0^{p+k}=v_0^{q+r}\) with the left-right
ordering are also as pictured in \tref{EndOfTrunk}.  Since the
numbered patterns for \(tv_0^k\) and \(\widehat{t}v_0^rs\) are
identical, the permutation given by \(s\) must be trivial on
\(\{0,1,\ldots, k\}\).  Thus there is a way to modify \(s\) using
the relations \tref{PiRelB}--\tref{PiRelD} from Lemma \ref{PiRels}
to give an \(s'\) that is a word in \(\{\sigma_i\mid i\ge k+1\}\).

Since there are \(p+k=q+r\) carets in the relevant trees,
there are \(p+k+1=q+r+1\) leaves numbered from \(0\) through
\(p+k=q+r\).  From Lemma \ref{PiToSigma} there is a way to modify
\(z'\) using the relations \tref{TwoVRelL}--\tref{TwoVRelQ} to give
a word \(z\) so that \(\Psi_{q+r}(z)=s'\).  Since there are no
appearances of \(\opi_i\) in \(z'\), there will be none in \(z\).

Now we use the fact that the only subscripts of the \(\sigma\) in
\(s'\) are above \(k\) to conclude that all the subscripts of the
\(\pi\) in \(z\) are below \((q+r-1)-k=(p+k-1)-k=p-1\).  Thus \(z\)
is a word in \(\{\pi_i\mid 0\le i\le p-2\}\).  This is what was
wanted.  

We now turn our attention to trees for \(L_1\).  We know that the
tree for \(\widehat{t}v_0^r\) below \(a'_m\) is pictured as in
\tref{EndOfTrunk}.  From Lemma \ref{RealizingIVOrder}, we know that
\(\widehat{t}v_0^r\) can be rewritten using the relations of \(\Pi\)
to end in \(k\) appearances of \(v_0\).  Thus we may assume that
\(j\ge k\).  In the word \(\widehat{t}v_0^rs'\) that corresponds to
\(L_1z\), we know that \(s'\) is a word in \(\{\sigma_i\mid i\ge
k+1\}\).  Thus the last \(k\) appearances of \(v_0\) in
\(\widehat{t}v_0^r\) can be moved to the right of \(s'\) using
relations \tref{PiRelE} at the expense of lowering each subscript in
\(s'\) by \(k\).  This changes \(s'\) to a word \(s''\) and changes
\(\widehat{t}v_0^rs'\) to \(\widehat{t}v_0^{r-k}s''v_0^k\).  Now we
have \[L_1z = (\widehat{t}v_0^rs, v_0^{q+r}) = (\widehat{t}v_0^rs',
v_0^{p+k}) = (\widehat{t}v_0^{r-k}s''v_0^k, v_0^{p+k}) =
(\widehat{t}v_0^{r-k}s'', v_0^p)\] and \(L_1\) can be represented by
\((\widehat{t}v_0^{r-k}, v_0^p)\).  Since \(\widehat{t}\) has \(q\)
letters, and \(p+k=q+r\), the word \(\widehat{t}v_0^{r-k}\) has
\(p\) letters. The trunk of the tree for \(\widehat{t}v_0^{r-k}\)
ends at vertex \(a'_m\) and so has \(m\) carets.  Thus setting
\(t'=\widehat{t}v_0^{r-k}\) completes the proof.  \end{proof}

The next lemma attacks non-normalized vertices in the trunk of a
tree.  We need a notion of complexity to measure progress.  In the
intermediate stages of the argument, it is extra work to define the
non-normalized vertices, so we use a different measure of progress.
The lemma will change the locations of labels, so we will focus on
the labels.  It turns out that only finitely many such changes can
be done, so this will be sufficient.

If \(T\) is a labeled tree, then we let \(a_0,a_1,\cdots, a_n\) be
the interior, left edge vertices of \(T\) reading from top to bottom
so that \(a_0\) is the root of \(T\).  We let \(b_0b_1\cdots b_n\)
be a word in \(\{0,1\}\) defined so that \(b_i=0\) if \(a_i\) is
labeled \(v\) and \(b_i=1\) if \(a_i\) is labeled \(h\).  We call
\(b_0b_1\cdots b_n\) the \emph{complexity} of \(T\).  If \(w_1\) and
\(w_2\) are two such words, then we say \(w_1<w_2\) if \(w_1\) is
shorter than \(w_2\) or if \(w_1\) and \(w_2\) are the same length
and \(w_1\) represents a smaller binary number than \(w_2\).  Note
that this gives the label of the root the most significant position.
The complexity will not be mentioned directly in the lemma, but will
be mentioned in its application.  However, the statement is easier
to understand if the complexity is kept in mind.

\begin{lemma}\mylabel{NormTrunk} let \(L=C_{i_0}C_{i_1}\cdots
C_{i_n}u\) where \(i_0<i_1<\cdots <i_n\) and \(u\) is a word of form
\(w(A,B)\).  Assume that the primary tree \(T\) for \(L\) is
normalized except at one or more vertices in the trunk of \(T\).
Let \(m\) be the number of carets in the trunk of \(T\).  Then
\(L\sim L'=C_{k_0}C_{k_1}\cdots C_{k_n}u'\) where
\(k_0<k_1<\cdots<k_n\), where \(u'\) is a word of form
\(w(A,B,\pi)\), so that \(m\ge k_n+1\), and so that the smallest
\(j\) so that \(i_j\ne k_j\) has \(i_j<k_j\).  \end{lemma}

\begin{proof}  Let \(\Lambda\) be the trunk of \(T\).  The interior
vertices of \(\Lambda\) are the interior, left edge vertices of
\(T\) and let these be \(a_0,a_1,\cdots, a_{m-1}\).  Let \(r\) be
the highest value with \(0\le r<m\) for which \(a_r\) is not
normalized.  Note that this is the lowest non-normalized interior
vertex of \(\Lambda\).

It follows that \(a_r\) has label \(h\) and that its children are
both interior vertices of \(T\) with label \(v\).  Because the left
child of \(a_r\) is an interior vertex, we must have \(r<m-1\).  The
vertex \(a_r\) must correspond to some \(C_{i_j}\) in \(L\) and from
Lemma \ref{LTreeStruct}, we have \(i_j=r\).

Represent \(T\) as \(\Lambda\) with a forest \(F\) attached.  Since
the left child of \(a_r\) is in \(\Lambda \) and has label \(v\), we
know that if \(j<n\), then \(i_{j+1}>r+1 = i_j+1\).  Since the right
child of \(a_r\) is an interior vertex, there is a letter in \(u\)
corresponding to it.  This letter must be some \(A_q\).  Also, the
right child of \(a_r\) is a root of \(F\).  By Lemma
\ref{RealizingIVOrder}, we can assume that \(A_q\) occurs as the
first letter of \(u\).  Thus we are looking at a word
\mymargin{WordInNormArg}\begin{equation}\label{WordInNormArg}
C_{i_0}C_{i_1}\cdots C_{i_j}C_{i_{j+1}}\cdots C_{i_n}A_q u''
\end{equation} where \(u''\) is the remainder of \(u\) after
\(A_q\).  The subword \(C_{i_0}C_{i_1}\cdots
C_{i_j}C_{i_{j+1}}\cdots C_{i_n}A_q\) of \tref{WordInNormArg} is a
trunk with a single caret labeled \(v\) attached at caret \(i_j\) of
the trunk on its right child.  From the details of the right-left
numbering, this implies that \(q=i_j\).  Thus \tref{WordInNormArg}
reads as
\mymargin{WordInNormArgII}\begin{equation}\label{WordInNormArgII}
C_{i_0}C_{i_1}\cdots C_{i_j}C_{i_{j+1}}\cdots C_{i_n}A_{i_j} u''.
\end{equation}

Since \(i_0<i_1<\cdots <i_n\) and \(i_{j+1}>i_j+1\), we can use
relations \tref{TwoVRelH} to rewrite \tref{WordInNormArgII} as 
\mymargin{WordInNormArgIII}\begin{equation}\label{WordInNormArgIII}
C_{i_0}C_{i_1}\cdots C_{i_j}A_{i_j}C_{i_{j+1}+1}\cdots
C_{i_n+1} u''.
\end{equation}
Combining relations \tref{TwoVRelI} and \tref{TwoVRelH} we get
\(C_mA_m=C_{m+1}B_m\pi_{m+1}\) so \tref{WordInNormArgIII} becomes 
\mymargin{WordInNormArgIV}\begin{equation}\label{WordInNormArgIV}
C_{i_0}C_{i_1}\cdots C_{i_j+1}B_{i_j}\pi_{i_j+1}C_{i_{j+1}+1}\cdots
C_{i_n+1} u''.
\end{equation}

Now \(i_{j+1}+1>i_j+2\), so relations \tref{TwoVRelJ} and
\tref{TwoVRelH} allow us to rewrite \tref{WordInNormArgIV} as 
\mymargin{WordInNormArgV}\begin{equation}\label{WordInNormArgV}
C_{i_0}C_{i_1}\cdots C_{i_j+1}C_{i_{j+1}}\cdots
C_{i_n}B_{i_j}\pi_{i_j+1} u''.
\end{equation}

If we now set \(B_{i_j}\pi_{i_j+1} u''= u'\) in
\tref{WordInNormArgV}, then \tref{WordInNormArgV} clearly satisfies
all provisions of the lemma except possibly \(m\ge k_n+1\).  However,
we observed above that \(i_j=r<m-1\) so \(i_j+2\le m\).  This gives
what we want in the case that \(j=n\).  If \(j<n\), then \(k_n=i_n\)
and \(m\ge i_n+1\) by Lemma \ref{LTreeStruct}.  \end{proof}

\begin{proof}[Proof of Proposition \ref{NormalLTree}] We start with
a word \(w\) in \(\Sigma_s\) of Section \ref{ConventionSec} and
assume \(w=LMR\) as specified in both  Lemmas \ref{ImprovMCor}  and
\ref{CsInFront}.  We work first on \(L\).  We have \(L\) expressible
as \((t, v_0^p)\) with \(t\) a word in \(w(v,h)\) of length \(p\)
and with the trunk of the tree \(T\) for \(t\) having \(m\) carets.

We apply Lemma \ref{ControlPerms} by letting \(L'=L\) in the
hypothesis of that lemma.  This gives \(L_1\) and \(L_2=L_1z\) with
\(L\sim L_2\), with \(z\) a word in \(\{\pi_i\mid i\le p-2\}\).
Also \(L_1\) is expressible as \((t', v_0^p)\) with \(t'\) a word of
form \(w(v,h)\) of length \(p\), with the trunk of the tree \(T'\)
for \(t'\) of length \(m\) and with \(T'\) normalized off the trunk.
Since we set \(L'=L\), we see that the trunks of \(T\) and \(T'\)
are identical.  Since the word \(z\) is in \(\{\pi_i\mid i\le
p-2\}\), it can be absorbed into \(M\) without disrupting our
assumptions on \(M\).  We now replace \(L\) with \(L_1\) and
proceed.

We now apply Lemma \ref{NormTrunk} to get \(L\sim L'\) as specified
in that lemma and then apply Lemma \ref{ControlPerms} to \(L\) and
\(L'\).  This gives \(L\sim L_2z\) as above.  The word \(z\) can be
added to \(M\).  The complexity of the tree \(T'\) as given in Lemma
\ref{ControlPerms} is dictated by the subword \(C_{k_0}C_{k_1}\cdots
C_{k_n}\) mentioned in that lemma.  But this is the sequence with
the same notation from Lemma \ref{NormTrunk} which gives a
complexity that is strictly less than that for the tree \(T\)
associated to \(L\).  Since there are finitely many complexities and
they are linearly ordered, the process of repeatedly applying Lemma
\ref{NormTrunk} followed by Lemma \ref{ControlPerms} must stop.  At
that point, the associated tree will be normalized.

To normalize the tree for \(R\), we apply this process to the
inverse of \(LMR\). \end{proof}

\subsection{An assumption of triviality}

We now look at our simplification of a word in the generators of
\(2V\) under the assumption that the word represents the trivial
element.  We will use the fact that this is also the trivial element
in \(\widehat{2V}\).

\begin{lemma}\mylabel{SamePatTrivPerm} Let \(w\) be a word in
\(\Sigma_s\) of Section \ref{ConventionSec} that represents the
trivial element of \(2V\).  Then \(w\sim 1\).  \end{lemma}

\begin{proof} Express \(w\sim LMR\) as in Lemma \ref{ImprovMCor} so
that \[LMR=(sv_0^{-p})(v_0^puv_0^{-p})(tv_0^{-P})^{-1}= sut^{-1}\]
where \(s\) and \(t\) are of form \(w(v,h)\) and \(u\) is of form
\(w(\sigma)\).  By Proposition \ref{NormalLTree}, we can assume that
the trees for \(t\) and \(s\) are normalized.  

Since \(sut^{-1} = (su,t)\) is the trivial element of
\(\widehat{2V}\), we know that the elements \(su\) and \(t\) are the
same elements of \(\Pi\) and represent the same numbered patterns.
Since \(s\) and \(t\) are words of form \(w(v,h)\) and \(u\) is a
word in \(\{\sigma_j\mid 0\le j\le p-1\}\), the words \(s\) and
\(t\) must give the same unnumbered pattern, and \(u\) must simply
renumber the numbered pattern for \(s\) to give the numbered pattern
\(su\).

The forests for \(s\) and \(t\) are trivial after the first trees and
so are normalized since the first trees are normalized.  Since the
forests for \(s\) and \(t\) are normalized and lead to the same
patterns, the forests are identical by Lemma \ref{UniqueNormFor}.
Since \(s\) and \(t\) are words of form \(w(v,h)\), the numbering on
the leaves is the left-right order and so \(s\) and \(t\) represent
the same numbered patterns.  Thus \(u\) affects the trivial
permutation on the numbering.  By Lemma \ref{ImprovMCor}, we know
\(M\sim 1\).

It remains to show that \(L\sim R^{-1}\).

The trunks of the trees for \(L\) and \(R\) are identical.  We have
\[L=C_{i_0}C_{i_1}\cdots C_{i_n}w(A,B)\] and
\[R^{-1}=C_{k_0}C_{k_1}\cdots C_{k_m}w'(A,B)\] with
\(i_0<i_1<\cdots<i_n\) and \(k_0<k_1<\cdots k_m\).  Since the
sequences \((i_0,i_1,\ldots, i_n)\) and \((k_0,k_1,\ldots k_m)\) are
determined by the labeling of the trunks of the trees for \(L\) and
\(R^{-1}\), they are identical as sequences.  What remains to be
shown is that \(w(A,B)\sim w'(A,B)\).  

The numbered, labeled forests \(F\) and \(F'\) for \(w(A,B)\) and
\(w'(A,B)\), respectively, that obtained by removing the trunks of
the trees for \(s\) and \(t\) are mirror images of the numbered,
labeled forests built from \(w(A,B)\) and \(w'(A,B)\) by replacing
each \(A\) by \(v\) and each \(B\) by \(h\).  Let these words be,
respectively, \(w(v,h)\) and \(w'(v,h)\).  Since the labeled forests
\(F\) and \(F'\) are equal, we know from Lemma \ref{TwoWdsSameFor}
that \(w(v,h)\) and \(w'(v,h)\) are related by relations
\tref{PiRelA}--\tref{PiRelE}.  However, these relations are
preserved under the homomorphism from \(\Pi\) to \(G\) under the
assignment \(v\mapsto A\) and \(h\mapsto B\).  Thus \(w(A,B)\sim
w'(A,B)\) and the proof is done.  \end{proof}

\subsection{The presentation}

Lemma \ref{SamePatTrivPerm} immediately gives the following.

\begin{thm}\mylabel{TwoVInfPres} The group \(2V\) is presented with
the generators of Section \ref{TwoVGenSec} and the relations
\tref{TwoVRelA}--\tref{TwoVRelQ}.  \end{thm}

\section{Finite presentations} 

We give finite presentations for \(\widehat{2V}\) and \(2V\).  We
proceed by showing that the relations that we have established for
\(\widehat{2V}\) and \(2V\) are consquences of finitely many of
those relations.  The techniques for doing this form a sort of
machine that would take about as long to describe as to use.  Thus
we do not make a theory out of it.  It is not clear that such a
theory is needed.

\subsection{A finite presentation for
\protect\(\widehat{2V}\protect\)} 

We have an infinite presentation from Theorem
\ref{WidehatTwoVPres}.  First we cut down the generating set
\(\{v_i, h_i, \sigma_i\mid i\in\N \}\).  The relations \tref{PiRelA}
and \tref{PiRelE} when \(i<j\) give \(v_i^{-1}x_jv_i=x_{j+1}\) when
\(x\) is any of \(v\), \(h\) or \(\sigma\).  This allows us to use
\[\begin{split} v_i &=v_0^{1-i}v_1v_0^{i-1}, \\ 
h_i &=v_0^{1-i}h_1v_0^{i-1}, \\ 
\sigma_i &=v_0^{1-i}\sigma_1v_0^{i-1}, \end{split}\] 
as definitions for all \(i\ge2\).  Thus \(\widehat{2V}\) is
generated by \(\{v_i, h_i, \sigma_i\mid i\in\{0,1\}\}\).

The relations \tref{PiRelA}--\tref{PiRelF} break into two classes:
the relations \tref{PiRelA}, \tref{PiRelC}, and \tref{PiRelE} whose
subscripts incorporate two parameters \(i\) and \(j\), and the
relations \tref{PiRelB}, \tref{PiRelD}, and \tref{PiRelF} whose
subscripts incorporate only the one parameter \(i\).  We treat these
two classes separately.

Conjugating the relations \tref{PiRelB}, \tref{PiRelD}, and
\tref{PiRelF} for \(i=1\) by powers of \(v_0\) shows that the
relations \tref{PiRelB}, \tref{PiRelD}, and \tref{PiRelF} for
\(i\ge2\) follow from \tref{PiRelB}, \tref{PiRelD}, and
\tref{PiRelF} for \(i=1\).  Thus \tref{PiRelB}, \tref{PiRelD}, and
\tref{PiRelF} are all consequences of the six relations obtained
when \(i\) is set to 0 or 1 in \tref{PiRelB}, \tref{PiRelD}, and
\tref{PiRelF}.

We consider \tref{PiRelA}.  If we rewrite \tref{PiRelA} when \(x\)
and \(y\) are both \(v\) as
\mymargin{PiRelAEquiv}\begin{equation}\label{PiRelAEquiv}
v_i^{-1}v_{i+k}v_i = v_{i+k+1}, \mathrm{\,\,for\,\, all\,\, }k>0,
\end{equation} then this is known to be true for \(i=0\) by
definition.  If \tref{PiRelAEquiv} is known for \(i=1\) and a set of
values of \(k>0\), then it is known for all \(i\) and that same set
of \(k\) values by conjugating the known expressions with \(i=1\) by
powers of \(v_0\).  If \tref{PiRelAEquiv} is known for \(i=1\) and
\(1\le k\le j\) with \(j\ge2\), then the calculation
\[v_1^{-1}v_{1+j+1}v_1 = v_1^{-1}v_{j+2}v_1 =
v_1^{-1}v_2^{-1}v_{j+1}v_2v_1 = v_3^{-1}v_{j+2}v_3 = v_{j+3}\] is
the inductive step that lets us conclude \tref{PiRelAEquiv} for
\(i=1\) and all \(k>0\) if we know \tref{PiRelAEquiv} for \(i=1\)
and \(k\in\{1,2\}\).  Thus \tref{PiRelAEquiv} for all needed \(i\)
and \(k\) follow from \[v_1^{-1}v_2v_1=v_3, \quad\mathrm{and}\quad
v_1^{-1}v_3v_1=v_4.\]

If we now look at \tref{PiRelA} when \(x=v\) and \(y=h\), then we
still have \(v_i^{-1}h_{i+k}v_i=v_{i+k+1}\) by definition when
\(i=0\).  Arguments similar to those in the previous paragraph show
that all of \tref{PiRelA} for \(x=v\) and \(y=h\) follow from
\[v_1^{-1}h_2v_1=h_3, \quad\mathrm{and}\quad v_1^{-1}h_3v_1=h_4.\]

When \(x=h\) in \tref{PiRelA}, we no longer have definitions to help
with \(i=0\).  Repetitions of the arguments above show that the
remaining cases of \tref{PiRelA} follow from the following eight
relations: \[\begin{split} h_0^{-1}v_1h_0 = v_2, &\qquad
h_1^{-1}v_2h_1=v_3, \\ h_0^{-1}v_2h_0=v_3,  &\qquad
h_1^{-1}v_3h_1=v_4, \end{split}
\qquad
\begin{split} h_0^{-1}h_1h_0 = h_2, &\qquad
h_1^{-1}h_2h_1=h_3, \\ h_0^{-1}h_2h_0=h_3, &\qquad
h_1^{-1}h_3h_1=h_4. \end{split}\]

The next set to consider is \tref{PiRelE}.  We give the discussion
and leave the results to be summarized after.  The relations
\tref{PiRelE} break into four smaller sets depending on the relative
sizes of \(i\) and \(j\).  For \(i=j\) and \(i=j+1\), we are dealing
with a single parameter, and we are in the same situation as
\tref{PiRelB}, \tref{PiRelD}, and \tref{PiRelF}.  Further, the fact
that the \(\sigma_i\) are their own inverses immediately shows that
the cases \(i=j\) and \(i=j+1\) give equivalent relations.  When
\(i<j\), we are in a situation ``isomorphic'' to some of the cases
in \tref{PiRelA}.  When \(i>j+1\), we are to prove that \(v_i\) or
\(h_i\) commutes with \(\sigma_j\).  When \(i\ge j+4\), the
calculation \[\sigma_j^{-1}v_i\sigma_j =
\sigma_j^{-1}v_{i-2}^{-1}v_{i-1}v_{i-2}\sigma_j =
v_{i-2}^{-1}v_{i-1}v_{i-2} = v_i\] gives the required inductive step.

Lastly, when \(i\ge j+4\), the inductive step
\[\sigma_j^{-1}\sigma_i\sigma_j =
\sigma_j^{-1}v_{i-2}^{-1}\sigma_{i-1}v_{i-2}\sigma_j = 
v_{i-2}^{-1}\sigma_{i-1}v_{i-2} = \sigma_i\] does the job for
\tref{PiRelC}.

Summarizing all of the above gives the following.  

\begin{thm}\mylabel{FiniteWHVPres}  The group \(\widehat{2V}\) is
presented by the generating set \(\{v_i, h_i, \sigma_i\mid
i\in\{0,1\}\}\) and the following set of 40 relations:
\begin{alignat*}{4}
v_2v_1&=v_1v_3,
  &v_3v_1&=v_1v_4,
  &h_2v_1&=v_1h_3,
  &h_3v_1&=v_1h_4, \\
v_1h_0&=h_0v_2, 
  &v_2h_0&=h_0v_3, 
  &v_2h_1&=h_1v_3, 
  &v_3h_1&=h_1v_4, \\
h_1h_0&=h_0h_2, 
  &h_2h_0&=h_0h_3, 
  &h_2h_1&=h_1h_3, 
  &h_3h_1&=h_1h_4, \\
\sigma_0v_2&=v_2\sigma_0,
  &\sigma_0v_3&=v_3\sigma_0,
  &\sigma_1v_3&=v_3\sigma_1,
  &\sigma_1v_4&=v_4\sigma_1, \\
\sigma_0h_2&=h_2\sigma_0,
  &\sigma_0h_3&=h_3\sigma_0, 
  &\sigma_1h_3&=h_3\sigma_1, 
  &\sigma_1h_4&=h_4\sigma_1, \\
\sigma_0v_0&=v_1\sigma_0\sigma_1,
  &\sigma_1v_1&=v_2\sigma_1\sigma_2,
  &\sigma_0h_0&=h_1\sigma_0\sigma_1,
  &\sigma_1h_1&=h_2\sigma_1\sigma_2, \\
\sigma_1h_0&=h_0\sigma_2, 
  &\sigma_2h_0&=h_0\sigma_3,
  &\sigma_2h_1&=h_1\sigma_3,
  &\sigma_3h_1&=h_1\sigma_4, \\
\sigma_2v_1&=v_1\sigma_3,
  &\sigma_3v_1&=v_1\sigma_4, 
  &\sigma_0^2&=1, 
  &\sigma_1^2&=1, \\
\sigma_0\sigma_2&=\sigma_2\sigma_0,
  &\sigma_0\sigma_3&=\sigma_3\sigma_0,
  &\sigma_1\sigma_3&=\sigma_3\sigma_1,
  &\sigma_1\sigma_4&=\sigma_4\sigma_1, \\
\sigma_0\sigma_1\sigma_0&=\sigma_1\sigma_0\sigma_1, \quad
  &\sigma_1\sigma_2\sigma_1&=\sigma_2\sigma_1\sigma_2, \quad
  &v_0h_1h_0 &=h_0v_1v_0\sigma1, \quad
  &v_1h_2h_1&=h_1v_2v_1\sigma_2.
\end{alignat*}
\end{thm}

\subsection{A finite presentation for
\protect\({2V}\protect\)} 

As in \(\widehat{2V}\) we get the following as definitions
\[\begin{split}
A_i&=A_0^{1-i}A_1A_0^{i-1}, \\
B_i&=A_0^{1-i}B_1A_0^{i-1}, \\
\pi_i&=A_0^{1-i}\pi_1A_0^{i-1}, \\
\opi_i&=A_0^{1-i}A\opi_1A_0^{i-1}, \\
C_m &=
(\opi_mB_m\opi_{m+1}\pi_m)(B_m\pi_{m+1}A_m^{-1})
\end{split}\] for \(i\ge2\) and \(m\ge0\).  The last line was
mentioned at the end of Section \ref{TwoVRelSec} and is easily
shown.  

Because of the homomorphism from \(\widehat{2V}\rightarrow 2V\)
determined by \(v\mapsto A\), \(h\mapsto B\), \(\sigma\mapsto\pi\),
all of the facts that we know about relations
\tref{PiRelA}--\tref{PiRelF} apply to the relations
\tref{TwoVRelA}--\tref{TwoVRelD}, \tref{TwoVRelK}-\tref{TwoVRelM},
and \tref{TwoVRelP} as mentioned in Lemma \ref{LMRABpi}.  Thus these
relations reduce to the image under the homomorhism
\(\widehat{2V}\rightarrow 2V\) of the relations in Theorem
\ref{FiniteWHVPres}.  If the remaining relations in the list
\tref{TwoVRelA}--\tref{TwoVRelQ} are treated in a manner similar to
that of Theorem \ref{FiniteWHVPres}, the following is proven.

\begin{thm}\mylabel{FiniteTVPres} The group \(2V\) is presented by
the generating set \(\{A_i, B_i, \pi_i, \opi_i\mid i\in\{0,1\}\}\)
the 40 relations obtained from the relations in Theorem
\ref{FiniteWHVPres} under the transformation \(v\mapsto A\),
\(h\mapsto B\), \(\sigma\mapsto \pi\) and the 30 relations below:
\begin{alignat*}{4}
\opi_2A_1&=A_1\opi_3, 
  &\opi_3A_1&=A_1\opi_4,
  &\opi_!B_0&=B_0\opi_2,
  &\opi_2B_0&=B_0\opi_3, \\
\opi_2B_1&=B_1\opi_3, 
  &\opi_3B_1&=B_1\opi_4,
  &\opi_0A_0&=\pi_0\opi_1,
  &\opi_1A_1&=\pi_1\opi_2, \\
\opi_0B_0&=C_1\pi_0\opi_1, \quad
  &\opi_1B_1&=C_2\pi_1\opi_2, \quad
  &C_2A_1&=A_1C_3, \quad
  &C_1A_1&=A_1C_4, \\
C_1B_0&=B_0C_2, 
  &C_2B_0&=B_0C_3,
  &C_2B_1&=B_1C_3,
  &C_3B_1&=B_1C_4, \\
C_0A_0&=B_0C_2\pi_1, 
  &C_1A_1&=B_1C_3\pi_2,
  &\opi_0^2&=1,
  &\opi_1^2&=1,  \\
\pi_0C_2&=C_2\pi_0,
  &\pi_0C_3&=C_3\pi_0,
  &\pi_1C_3&=C_3\pi_1,
  &\pi_1C_4&=C_4\pi_1, \\
\pi_0\opi_2&=\opi_2\pi_0,
  &\pi_0\opi_3&=\opi_3\pi_0,
  &\pi_1\opi_3&=\opi_3\pi_1,
  &\pi_1\opi_4&=\opi_4\pi_1, \\
\pi_0\opi_1\pi_0&=\opi_1\pi_0\opi_1,
  &\pi_1\opi_2\pi_1&=\opi_2\pi_1\opi_2.
\end{alignat*}
\end{thm}

\section{Normal forms}

The group \(2V\) and \(V\) resemble each other greatly, but differ
in one important aspect.  Both \(2V\) and \(V\) live in larger
groups of fractions, respectively \(\widehat{2V}\) and
\(\widehat{V}\), which in some ways are better behaved than \(2V\)
and \(V\).  However, \(\widehat{V}\) has a property not possessed by
\(\widehat{2V}\).  

Elements of \(\widehat{V}\) have a nice normal form when regarded as
pairs.  The positive monoid \(P\) for \(\widehat{V}\) has a nice
length function, and a given element of \(\widehat{V}\) has a nice
class of pairs \((Y,Z)\) that represent it that is characterized by
the fact that they are minimal in length in the monoid \(P\).  In
particular, if \((Y,Z)\) and \((Y',Z')\) are in this class, then
\(Y\) and \(Y'\) (and also \(Z\) and \(Z'\)) differ by an invertible
element (a permutation) in \(P\).  This gives \(\widehat{V}\) a nice
semi-normal form that can be turned into a normal form by an easy
shift of point of view.  See Theorem 2 in Section 6.3 of
\cite{brin:bv}.

The remarks above about \(\widehat{V}\) follow from a property of
\(P\) called least common left multiples in \cite{brin:bv}.  A
common left multiple of two elements \(Y\) and \(Z\) in \(P\) is an
element \(L=AY=BZ\) with \(A\) and \(B\) in \(P\).  If \(L\) and
\(L'\) are two left multiples of \(Y\) and \(Z\), then we write
\(L\le L'\) if there is a \(C\) in \(P\) so that \(L'=CL\).  A least
common left multiple of \(Y\) and \(Z\) is a common left multiple of
\(Y\) and \(Z\) that is least in this order.  A monoid had
\emph{least common left multiples} if every pair in the monoid with
a common left multiple has a least common left multiple.  The
definition is worded so that a monoid can have least common left
multiples even if not every pair has a common left multiple.

The semi-normal form discussed above for elements of
\(\widehat{V}\) comes directly from the fact that the positive
monoid \(P\) of \(\widehat{V}\) has least common left multiples.
The monoid \(\Pi\) of this paper does not.  The culprit is relation
\tref{PiRelF}.  We leave it as an exercise to show that while
\(v_0\) and \(h_0\sigma_1\) have common left multiples
\[(h_0v_1\sigma_2)v_0 = h_0v_1v_0\sigma_3 =
v_0h_1h_0\sigma_1\sigma_3 = (v_0h_1\sigma_2)(h_0\sigma_1)\] and
\[(h_0v_1)v_0 = (v_0h_1)(h_0\sigma_1),\] there is no least common
left multiple for \(v_0\) and \(h_0\sigma_1\).  From this it follows
that the element \[(h_0v_1\sigma_2,\,\, h_0v_1) = (v_0h_1\sigma_2,
\,\, v_0h_1)\] of \(\widehat{2V}\) has no unique ``minimal''
representative.  This accounts for the arbitrary choice made (a
normalized vertex with a secondary label must have label \(v\)) in
the definition of a normalized forest.

It would be nice to know if \(\widehat{2V}\) and \(2V\) are of type
\(F_\infty\) (have classifying spaces that are finite in each
dimension).  It seems at the moment that the absence of the least
common left multiples property will make the question harder.

\vspace{20pt}


\begin{thebibliography}{10}

\bibitem{boardman+vogt}
J.~M. Boardman and R.~M. Vogt, \emph{Homotopy invariant algebraic structures on
  topological spaces}, Springer-Verlag, Berlin, 1973, Lecture Notes in
  Mathematics, Vol. 347. \MR{54 \#8623a}

\bibitem{brin:hd3}
Matthew~G. Brin, \emph{Higher dimensional {T}hompson groups}, Geometriae
  Dedicata, to appear.

\bibitem{brin:zs}
\bysame, \emph{On the {Z}appa-{S}z\'ep product}, Communications in Algebra, to
  appear.

\bibitem{brin:bv}
\bysame, \emph{The algebra of strand splitting. {I}. {T}he algebraic structure
  of the braided {T}hompson group.}, preprint, Binghamton University, 2004.

\bibitem{brin:bv3}
\bysame, \emph{The algebra of strand splitting. {II}. {A} presentation for the
  braid group on one strand}, preprint, Binghamton University, 2004.

\bibitem{CFP}
J.~W. Cannon, W.~J. Floyd, and W.~R. Parry, \emph{Introductory notes on
  {R}ichard {T}hompson's groups}, Enseign. Math. (2) \textbf{42} (1996),
  no.~3-4, 215--256. \MR{98g:20058}

\bibitem{cliff+prest:I}
A.~H. Clifford and G.~B. Preston, \emph{The algebraic theory of semigroups.
  {V}ol. {I}}, American Mathematical Society, Providence, R.I., 1961,
  Mathematical Surveys, No. 7.

\bibitem{guba+sapir2}
Victor Guba and Mark Sapir, \emph{Diagram groups}, Mem. Amer. Math. Soc.
  \textbf{130} (1997), no.~620, viii+117. \MR{98f:20013}

\bibitem{may:geom+loop}
J.~P. May, \emph{The geometry of iterated loop spaces}, Springer-Verlag,
  Berlin, 1972, Lectures Notes in Mathematics, Vol. 271. \MR{54 \#8623b}

\bibitem{squier:word-problems}
Craig~C. Squier, \emph{Word problems and a homological finiteness condition for
  monoids}, J. Pure Appl. Algebra \textbf{49} (1987), no.~1-2, 201--217.
  \MR{89a:20059}

\end{thebibliography}

\providecommand{\bysame}{\leavevmode\hbox to3em{\hrulefill}\thinspace}
\providecommand{\MR}{\relax\ifhmode\unskip\space\fi MR }
\providecommand{\MRhref}[2]{%
  \href{http://www.ams.org/mathscinet-getitem?mr=#1}{#2}
}
\providecommand{\href}[2]{#2}

\noindent Department of Mathematical Sciences

\noindent State University of New York at Binghamton

\noindent Binghamton, NY 13902-6000

\noindent USA

\noindent email: matt@math.binghamton.edu

\end{document}